\documentclass[12pt,a4paper]{article}
\usepackage{latexsym}
\usepackage{graphics}
\newtheorem{theorem}{Theorem}
\newtheorem{alg}{Algorithm}
\newtheorem{lemma}{Lemma}
\newtheorem{cor}{Corollary}
\newcounter{step}
\newcommand{\balg}
 {\begin{list}{Step \thestep.\hfill}
  {\usecounter{step}
   \setlength{\leftmargin}{4em}
   \setlength{\labelwidth}{\leftmargin}
   \addtolength{\labelwidth}{-1.0\labelsep}
   \setlength{\rightmargin}{\leftmargin}      }    }
\newcommand{\ealg}{\end{list}}
\newcommand{\beq}{\begin{equation}}
\newcommand{\eeq}{\end{equation}}

\title{Piecewise Convex-Concave Approximation in the $\ell_{\infty}$ Norm}
\author{M. P. CULLINAN  \\ Maryvale Institute, Birmingham, B44 9AG, UK
\\ mpcullinan@gmail.com
}
\date{26 July 2010}
\newcommand{\reals}{{\hbox{I\kern-.2em\hbox{R}}}}
\newcommand{\union}{\cup}
\newcommand{\ff}{\mbox{$\bf f$}}
\newcommand{\realsn}{\mbox{$\reals^{n}$}}
\newcommand{\linf}{\mbox{$\ell_{\infty}$}}
\newcommand{\Yq}{\ensuremath{Y_{q}}}
\newcommand{\vv}{\mbox{$\bf v$}}
\newcommand{\civ}[2]{\mbox{$c_{#1}(\bf #2)$}}
\newcommand{\yy}{\mbox{$\bf y$}}
\newcommand{\hq}{\mbox{$h_{q}$}}
\newcommand{\iplus}[2]{\ensuremath{I^{+} ({#1 , #2})}}
\newcommand{\iminus}[2]{\ensuremath{I^{-} ({#1 , #2})}}
\newcommand{\isig}[2]{\ensuremath{I^{\sigma} ({#1 , #2})}}

\newcommand{\beqa}{\begin{eqnarray}}
\newcommand{\eeqa}{\end{eqnarray}}
\newcommand{\beqax}{\begin{eqnarray*}}
\newcommand{\eeqax}{\end{eqnarray*}}
\newcommand{\hsig}[2]{\ensuremath{h^{\sigma} ({#1 , #2})}}
\newcommand{\oneton}[1]{\mbox{$1 \leq #1 \leq n$}}
\newcommand{\half}{\mbox{$\frac{1}{2}$}}
\newcommand{\yzero}{\ensuremath{{\bf y}^{0}}}
\newcommand{\hzero}{\ensuremath{h_{0}}}

\newcommand{\hb}{\mbox{$\overline{h}$}}
\newcommand{\gb}{\ensuremath{\overline{g}}}
\newcommand{\Ib}{\ensuremath{\overline{I}}}
\newcommand{\Sb}{\ensuremath{\overline{S}}}
\newcommand{\doblob}[2]{\put(#1,#2){\circle*{4}}}
\newcommand{\dosmblob}[2]{\put(#1,#2){\circle*{2}}}
\newcommand{\writeat}[3]{\put(#1,#2){\makebox(0,0){#3}}}
\newcommand{\doline}[4]{\put(#1,#2){\line(#3){#4}}}
\newcommand{\yone}{\ensuremath{{\bf y}^{1}}}
\newcommand{\bea}{\begin{eqnarray*}}
\newcommand{\eea}{\end{eqnarray*}}

\newcommand{\yyb}{\mbox{$\overline{\yy}$}}
\newcommand{\yb}{\ensuremath{\overline{y}}}
\newcommand{\stoalpha}{\mbox{$\{ \, s,s_{1}, \ldots , s_{\alpha} \, \}$}}

\newcommand{\jst}{\ensuremath{j^{\star}}}
\newcommand{\fjs}{\mbox{$f_{\scriptstyle j^{\star}}$}}
\newcommand{\kpl}{\ensuremath{k^{+}}}
\newcommand{\sbkp}[1]{\mbox{$#1_{k^{+}}$}}
\newcommand{\rfeqs}[2]{(\ref{#1})--(\ref{#2})}

\newcommand{\alphab}{\ensuremath{\overline{\alpha}}}

\newcommand{\mam}{\mbox{$(-1)^{\alpha -1}$}}

\newcommand{\qb}{\ensuremath{\overline{q}}}

\newcommand{\Pal}{\mbox{$P_{\alpha}$}}

\newcommand{\sal}{\ensuremath{s_{\alpha}}}
\newcommand{\salp}{\ensuremath{s'_{\alpha}}}
\newcommand{\talm}{\ensuremath{t_{\alpha - 1}}}
\newcommand{\talpm}{\ensuremath{t'_{\alpha - 1}}}
\newcommand{\tal}{\ensuremath{t_{\alpha}}}
\newcommand{\talp}{\ensuremath{t'_{\alpha}}}
\newcommand{\salm}{\ensuremath{s_{\alpha - 1}}}

\newcommand{\talpl}{\mbox{$t_{\alpha+1}$}}

\newcommand{\csts}[1]{\mam c_{s(\alpha),t(\alpha),s(\alpha +1)}(#1)}
\newcommand{\ctst}[1]{\mam c_{t(\alpha),s(\alpha +1),t(\alpha +1)}(#1)}
\newcommand{\cstt}[1]{\mam c_{s(\alpha),t(\alpha),t(\alpha +1)}(#1)}
\newcommand{\csst}[1]{\mam c_{s(\alpha),s(\alpha+1),t(\alpha +1)}(#1)}
\newcommand{\ggg}{\mbox{$\bf g$}}
\newcommand{\Pu}{\mbox{$P_{\infty}$}}
\newcommand{\Pe}{\mbox{$P_{2}$}}
\newcommand{\eps}{\mbox{$\epsilon$}}
\begin{document}
\maketitle
\begin{abstract}
Suppose that $\ff \in \reals^{n}$ is a vector of $n$ error-contaminated measurements of $n$ smooth values measured at distinct and strictly ascending abscissae. The following projective technique is proposed for obtaining a vector of smooth approximations to these values. Find \yy\ minimizing $\| \yy - \ff \|_{\infty}$ subject to the constraints that the second order consecutive divided differences of the components of \yy\ change sign at most $q$ times. This optimization problem (which is also of general geometrical interest) does not suffer from the disadvantage of the existence of purely local minima and allows a solution to be constructed in $O(nq)$ operations. A new algorithm for doing this is developed and its effectiveness is proved. Some of the results of applying it to undulating and peaky data are presented, showing that it is economical and can give very good results, particularly for large densely-packed data, even when the errors are quite large.
\end{abstract}
%
%
\section{Introduction}
Given a set of observations $f_{j}$,
at strictly ascending absciss\ae\ $x_{j}$, for $1 \leq j \leq n$,
it may be known that the
observations represent measurements of smooth quantities
contaminated by errors.  A method is then needed for obtaining a
smooth set of points while respecting the observations as much as possible.
One method is to make the least change to the observations,
measured by a suitable norm,
in order to achieve a prescribed definition of smoothness.
The data smoothing method of Cullinan \& Powell (1982)
proposes defining smoothness as the consecutive divided differences
of the points of a prescribed order $r$ having at most a
prescribed number $q$ of sign changes.
There are many good reasons for using the number of sign changes in the divided differences of data as a criterion of smoothness. Normally data values of a smooth function will have very few sign changes, whereas if even one error is introduced, it will typically cause $k$ sign changes in the $k$th order divided differences of the contaminated data. Thus constructing a table of divided differences is a cheap and sensitive test for smoothness (see, for example, Hildebrand, 1956).
If the observations $f_{j}$, for $1 \leq j \leq n$, are regarded as the
components of a vector $\ff \in \reals^{n}$ and the function
$F : \reals^{n} \rightarrow \reals$ is defined through the chosen norm by
$F( {\bf v} ) = \| {\bf v - f} \|$, then the data smoothing problem becomes the constrained
minimization of $F$.
This approach has several advantages. There is no need to choose (more or less arbitrarily) a set of approximating functions, indeed the data are treated as the set of finite points which they are rather than as coming from any underlying function. The method is projective or invariant in the sense that it leaves smoothed points unaltered. It depends on two integer parameters which will usually take only a small range of possible values, rather than requiring too arbitrary a choice of parameters.  It may be possible to choose likely values of $q$ and $r$ by inspection of the data. The choice of norm can sometimes be suggested by the kind of errors expected, if this is known. For example the $\ell_{1}$  norm is a good choice if a few very large errors are expected, whereas the \linf\ norm might be expected to deal well with a large number of small errors. There is also the possibility that the algorithms to implement the method may be very fast.

The main difficulty in implementing this method is that when
$q \geq 1$, the possible existence of purely local minima of $F$ makes
the construction of an efficient algorithm very difficult.
This has been done for the $\ell_{2}$ norm for $r \leq 2$---see for example Demetriou (2002). The author dealt with the case $q=0$ and arbitrary $r$ for the $\ell_{2}$ norm (Cullinan, 1990).

It was claimed in Cullinan \& Powell (1982) that when the \linf\ norm is
chosen and $r=2$, all the local minima of $F$ are global
and a best approximation can be constructed in $O(nq)$ operations;
and an algorithm  for doing this was outlined.
These claims were proved by Cullinan (1986) which also considered the case $r=1$. It was
shown that in these cases the minimum value of $F$ is determined by $q+r+1$ of the
data, and a modified algorithm for the case $r=2$ was developed which is believed to
be better than that outlined in Cullinan \& Powell (1982). This
new algorithm will now be presented in Section 2 and its effectiveness
will be proved. Section 3 will then describe the results of some
tests of this method which show that it is a very cheap way of
filtering noise but can be prone to end errors.
%
%
%
\section{The algorithm}
 This section will construct a best \linf\ approximation to a vector $\ff \in \realsn$ with not
more than $q$ sign changes in its second divided differences. More
precisely let $\vv \in \realsn$ with $x_{1} < x_{2} < \ldots < x_{n}$ and
 \beq F(\bf{v}) = \| {\bf v-f}\|_{\infty}
 \label{210}
 \eeq
 \beq
 c_{ijk} ({\bf v}) =
 \frac{1}{x_{k}-x_{i}}
 \left(
  \frac{v_{k}-v_{j}}{x_{k}-x_{j}} - \frac{v_{j}-v_{i}}{x_{j}-x_{i}}
  \right),
\label{211}
\eeq
\beq \civ{i}{v} = c_{i,i+1,i+2}({\bf v}). \label{213} \eeq

The feasible set of points $\Yq \subset \realsn$ is defined as the
set of all vectors $\vv \in \realsn$ for which the signs of the
successive elements of the sequence $1,\civ{1}{v}, \ldots
,\civ{n-2}{v}$ change at most $q$ times, and the the problem is
then to develop an algorithm to minimize $F$ over $\Yq$.

\subparagraph{}

The solution depends on the fact that the value \hq\ of the best
approximation is determined by $q+3$ of the data. Since a best
\linf\ approximation is not unique, there is some choice of which
one to construct. The one chosen, $\yy$, has the following
property: $y_{1}=f_{1}$, $y_{n}=f_{n}$, and for any $j$ with $2
\leq j \leq n-1$,
 \[
 \mbox{if } \pm \civ{j-1}{y} > 0 \mbox{  then  } y_{j}=f_{j} \pm \hq\
 \]
The vector \yy\ is then determined from $\hq$, from the set of
indices $i$ where $\civ{i-1}{ \yy} \neq 0$, and from the ranges
where the divided differences do not change sign.

\subparagraph{}

The method by which the best approximation is constructed and the
proof of the effectiveness of the algorithm that constructs it
are best understood by first considering the cases $q=0,1,$ and
$2$ in detail and giving algorithms for the construction of a
best approximation in each case. Once this has been done it is
easy to understand  the general algorithm.

When $q=0$, this best approximation is formed from the ordinates
of the points on the lower part of the boundary of the convex
hull of the points $(x_{j},f_{j})$, for $1 \leq j \leq n$, (the
graph of the data in the plane) by increasing these ordinates by
an amount $h$.

When $q=1$, there exist integers $s$ and $t$ such that $1 \leq s
\leq t \leq n$, and $y_{1}, \ldots ,y_{s}$ are the ordinates on
the lower part of the boundary of the convex hull of the data
$f_{1}, \ldots , f_{s}$ increased an amount $h$; $y_{t}, \ldots
,y_{n}$ are ordinates on the upper part of the boundary of the
convex hull of the data $f_{t}, \ldots , f_{n}$ decreased an
amount $h$; and if $s < j < t$, $y_{j}$ lies on the straight line
joining $(x_{s}, y_{s})$ to $(x_{t}, y_{t})$. The best
approximation therefore consists of a convex piece and a concave
piece joined where necessary by a straight line.

The best approximation over \Yq\ consists of $q+1$ alternately
raised pieces of lower boundaries of convex hulls and lowered
pieces of upper boundaries of concave hulls joined where necessary
by straight lines. These pieces are built up recursively from
those of the best approximation over $Y_{q-2}$.

The points on the upper or lower part of the boundary of the
convex hull of the graph of a range of the data each lie on a
convex polygon and are determined from its vertices. The
algorithms to be described construct sets of the indices of these
vertices and the value of a best approximation. The best
approximation vector is then constructed by linear interpolation.

\subparagraph{}

Before considering the cases $q=0$, $1$, and $2$, an important
preliminary result will be established. It is a tool that helps
to show that the vectors constructed by the algorithms are
optimal. The value \hq\ of the best approximation over \Yq\ will
be found by the algorithm, together with a vector $\yy \in\Yq$
such that $F(\yy) = \hq$. To show that \yy\ is optimal, a set $K$
of $q+3$ indices will be constructed such that if \vv\ is any
vector in \realsn\ such that $F(\vv) < \hq$, then the consecutive
second divided differences of the components $v_{k}$, for $k \in
K$, change sign $q$ times starting with a negative one. It will
then be inferred that $\vv \not\in \Yq$. In order to make this
inference it must be shown that the consecutive divided
differences of {\em all} the components of \vv\ have at least as
many sign changes as those of the components with indices in $K$.
This result will now be proved.
\begin{theorem}
Let $K \subseteq \{ \, 1, \ldots, n \, \}$ and let $\vv \in
\realsn$ be any vector such that the second divided differences of
the $v_{k}$, for $k \in K$, change sign $q$ times. Then the
divided differences of all the components of \vv\ change sign at
least $q$ times.
\end{theorem}
\paragraph{Proof}
Firstly, suppose that $K$ is formed by deleting one element $j$
from $\{ \, 1, \ldots, n \, \}$ and that $3 \leq j \leq n-2$. Let
$c_{j-2}, c_{j-1}, c_{j}$ be defined from \vv\ by (\ref{211}) and
(\ref{213}), and let $c^{\prime}_{j-2} = c_{j-2, j-1, j+1}(\vv)$
and $c^{\prime}_{j-1} = c_{j-1, j+1, j+2}(\vv)$ denote the new
divided differences that result from deleting $j$. Manipulation
of (\ref{211}) yields the equations
\[
c^{\prime}_{j-2} =
    \frac{x_{j}- x_{j-2}}{x_{j+1} - x_{j-2}} c_{j-2} +
    \frac{x_{j+1}- x_{j}}{x_{j+1} - x_{j-2}} c_{j-1}
\]
and
\[
c^{\prime}_{j-1} =
    \frac{x_{j+2}- x_{j}}{x_{j+2} - x_{j-1}} c_{j} +
    \frac{x_{j}- x_{j-1}}{x_{j+2} - x_{j-1}} c_{j-1} ,
        \]

\noindent so that $c^{\prime}_{j-2}$ lies between $c_{j-2}$ and
$c_{j-1}$ and $c^{\prime}_{j-1}$ lies between $c_{j-1}$ and
$c_{j}$. It follows that the number of sign changes in the
sequence \beq \ldots, c_{j-3}, c_{j-2}, c^{\prime}_{j-2}, c_{j-1},
  c^{\prime}_{j-1}, c_{j}, c_{j+1}, \ldots ,
\label{214} \eeq is the same as that in the sequence \[\ldots,
c_{j-3}, c_{j-2},  c_{j-1}, c_{j}, c_{j+1} \ldots , \] and hence
that deleting $c_{j-2},  c_{j-1}, c_{j}$ from (\ref{214}) cannot
increase the number of sign changes. The same argument covers the
cases $j=2$ and $j=n-1$, and when $j=1$ or $j=n$ this result is
immediate.

Repeated application of this result as elements $j$ of the set
$\{ \, 1, \ldots, n \, \}  \backslash K$ are successively deleted
from $ \{ \, 1, \ldots, n \, \}$ then proves the theorem.
~$\quad\Box$

\begin{cor}
If $i \leq k-2$ and all the divided differences $c_{j}(\vv), i <
j < k$, are non-negative (or non-positive), and if $i \leq r<s<t
\leq k$, then $c_{rst}(\vv)$ is also non-negative (or
non-positive).
\end{cor}

This Theorem is crucial to the effectiveness of the algorithm. In
Cullinan (1986) it was proved that the set of best approximations
is connected, so that purely local minima are ruled out, but that
this is not the case for higher orders of divided differences.
This Theorem allows the explicit construction of global minima
determined by $q+3$ of the data (and so it does not seem
necessary to prove connectedness). There is no analogous result
for higher order divided differences, and so no ready
generalization of the methods of this paper to such cases.
%
%
%
\subsection{The case $q=0$}
When $q=0$, the required solution is a best convex approximation
to \ff. The particular one, $\yzero$, that will be constructed
here was first produced by Ubhaya (1979). It will also be
convenient to construct a best \emph{concave} approximation to
data.

The convex approximation is determined from the vertices of the
lower part of the boundary of the convex hull of the points, and
the concave approximation from the vertices of the upper part.
Each of these sets of indices can be constructed in $O(n)$
operations by the following algorithm. It is convenient to apply
the construction to a general range of the data and to describe
it in terms of sets of indices.
Accordingly, define a \emph{range}
    $[ r , s ] = \{ j : r \leq j \leq  s \}$,
and a \emph{vertex set} of $[r,s]$ to be any set $I$ such that
    $\{ r,s\} \subseteq I \subseteq [r,s]$.
Given a vertex set $I$ of a range $[r,s]$ and also quantities
$v_{i}$, $i \in I$, define the \emph{gradients}
\[
    g_{ik}({\bf v}) =
                    \frac{v_{k} - v_{i}}{x_{k} - x_{i}}
                   \qquad \mbox{ for }
                    i \neq k .
\]

Given any index $j : 1 \leq j \leq n$, define the
\emph{neighbours} of $j$ in $I$ by
    \beqax
    j^{+} (I) & = & \min_{i \in I} \{i>j \} \mbox{ when } j<s
    \\
    j^{-} (I) & = & \max_{i \in I} \{ i<j \} \mbox{ when } j>r,
    \eeqax
and (for purposes of extrapolation)
    \beqax
    j^{+} (I) & = & s^{-}(I)\mbox{ when } j \geq s
    \\
    j^{-} (I) & = & r^{+} (I)\mbox{ when } j \leq r .
    \eeqax
The \emph{interpolant} $\vv (I)$ can now be defined by
    \beqa
    v_{j} (I) & = & v_{j} \qquad \mbox{for } j \in I
    \label{2103a}
    \\
    v_{j} (I) & = & v_{j^{-}(I)} + g_{j^{-}(I) \, j^{+}(I)} ({\bf f})
                                (x_{j} - x_{j^{-}(I)} )\qquad \mbox{for }j \notin I.
    \label{2103b}
    \eeqa
When $I= \{ p,q \}$ it is convenient to write $v_{j} (I)$ as
$v_{j} (p,q)$ etc.

The two cases of convex and concave approximations are handled
using the sign variable $\sigma$, where $\sigma = +$  for the
convex case and $\sigma = -$  for the concave case. The
\emph{convex} and \emph{concave} \emph{optimal vertex sets}
\iplus{r}{s} and \iminus{r}{s} are then constructed by systematic
deletion as follows.

\begin{alg} To find \isig{r}{s} when $r \leq s-2$.
\end{alg}
\balg
\item  Set $I:= [r,s]$ and $i:=r$, $j:=r+1$, $k:=r+2$.
\item  Evaluate $c := c_{ijk}(\ff)$.  If $\sigma c > 0$: go to Step~5.
\item  Delete $j$ from $I$. If $i=s$: go to Step~5.
\item  Set $k := j$, $j := i$, and $i := i^{-}(I)$. Go to Step~2.
\item  If $k=t$: set $\isig{r}{s}=I$ and {\bf stop.} \\
       Otherwise: set $i := j$, $j := k$, and $k := k^{+}(I)$. Go to Step~2.
\ealg

The \emph{price} of making a convex/concave approximation in the
range $[ r,s ]$ is given by \beq
    \hsig{r}{s} = \half \max_{j \in [ r,s ]} \sigma
                    ( f_{j} - f_{j}(\isig{r}{s}),
\label{2101}
 \eeq
and the required best approximation \yzero\ is then defined by
\beq
    y_{j}^{0} = f_{j} (I) + h,
    \label{2102}
    \eeq
where $I=I^{+}(1,n)$ and $h=h^{+}(1,n)$.

The construction of \yzero\ is illustrated in Figure~1.
%
%
Fig 101
\begin{figure}
\begin{center}
\begin{picture}(334,196)
\thicklines \put(0,180){\circle*{4}} \put(0,180){\line(1,-2){36}}
\put(36,108){\circle*{4}}
\put(45,120){\makebox(0,0){$f_{k^{-}}(I)$}}
\put(36,96){\makebox(0,0){$k^{-}$}} \put(36,108){\line(1,-1){72}}
\put(108,36){\circle*{4}} \put(108,48){\makebox(0,0){$f_{k}(I)$}}
\put(108,24){\makebox(0,0){$k$}} \put(108,36){\line(1,0){108}}
\put(180,36){\circle*{2}} \put(180,48){\makebox(0,0){$f_{j}(I)$}}
\put(180,24){\makebox(0,0){$j$}} \put(180,120){\circle*{4}}
\put(180,132){\makebox(0,0){$f_{j}$}} \put(216,36){\circle*{4}}
\put(216,48){\makebox(0,0){$f_{k^{+}}(I)$}}
\put(216,24){\makebox(0,0){$k^{+}$}} \put(216,36){\line(3,2){54}}
\put(270,72){\circle*{4}} \put(270,72){\line(1,2){36}}
\put(306,144){\circle*{4}}
\end{picture}
\end{center}
\caption{Consecutive elements $k^{-}$, $k$, $k^{+}$ of
$I=I^{+}(1,n)$ and the construction of $f_{j}(I)$}
\end{figure}
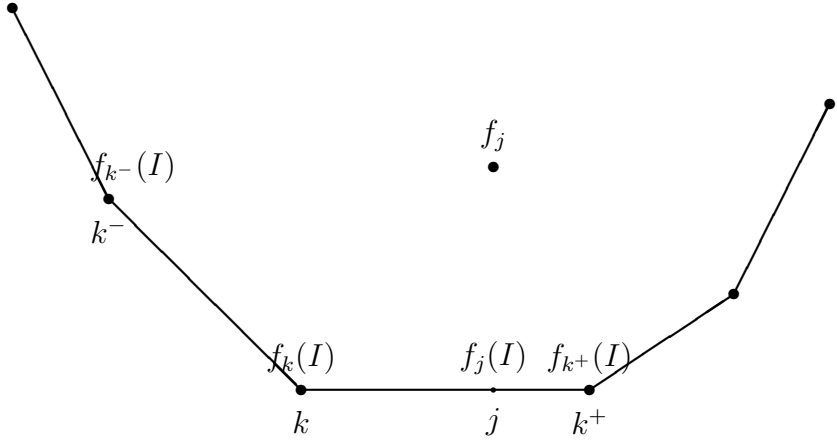

%

\begin{theorem}
Let $n \geq 3$, and let $\hzero= h^{+}(1,n)$ and \yzero\ be
defined from (\ref{2101}) and (\ref{2102}) using Algorithm 1,
(\ref{2103a}) and (\ref{2103b}). Then
    $F( \yzero )= \hzero = \inf F(Y_{0})$
and $\yzero \in Y_{0}$.
\end{theorem}
\paragraph{Proof}
The first remark is that Algorithm 1 produces a well-defined
vertex set $I$ of $[r,s]$ from which the quantities $f_{j}(I)$
are also well-defined for all $j$, so that $\hzero$ and \yzero\
are also well-defined.

The proof that the points $(x_{j}, f_{j}(I))$, \oneton{j}, lie on
the lower part of the convex hull of the data
is in Ubhaya. The components $f_{j}(I)$ can be defined as those
that are maximal subject to the inequalities
    \beq
    f_{j}(I) \leq f_{j}, \qquad
        \mbox{for }\oneton{j}, \label{2105}
    \eeq
and
    \beq
    c_{i}({\bf f}(I)) \geq 0, \qquad
        \mbox{for } 1 \leq i \leq  n-2, \label{2106}
    \eeq
{\em i.e.}\ ${\bf f}(I)\in Y_{0}$. It follows that
    $\yzero \in Y_{0}$
 and from (\ref{2101}) and (\ref{2105}) that $F( \yzero )= h$.

It remains to prove that \yzero\ is optimal. If $h =0$, \yzero\
must be optimal. If $h>0$, there will be a lowest integer
    $\jst \not\in I$
such that equality is attained in (\ref{2101}), and since $1$ and
$n$ can never be deleted from $I$, it must be the case that $2
\leq \jst \leq n-1$. Then
\[
y^{0}_{\jst } = f_{\jst} - \hzero ,
\]
and $k=j^{-}(I)$ and $k^{+}= j^{+}(I)$ are consecutive
elements of $I$ such that
    $k<j < k^{+}$
and
\[
y^{0}_{k} = f_{k} + \hzero ,
\]
\[
y^{0}_{k^{+}} = f_{k^{+}} + \hzero .
\]
Since $\jst \not\in I$, $c_{k, \jst ,k^{+}}( \yzero\ )
    = c_{k, \jst ,k^{+}} ( {\bf f}(I)) =0$.
Now if $\vv \in \realsn$ and $F(\vv) < h$, then
    \[
    v_{k} < f_{k} + h =y^{0}_{k},
    \]
    \[
    v_{\jst } > f_{\jst } - h = y^{0}_{ \jst },
    \]
    \[
    v_{k^{+}} < f_{k^{+}} + h =y^{0}_{k^{+}}.
    \]

Now the constraint function $c_{k, \jst, k^{+}}$
is an increasing function of $v_{k}$ and $v_{k^{+}}$  and a
decreasing function of $v_{\jst }$, so
    $c_{k, \jst ,k^{+}}(\vv ) < c_{k, \jst ,k^{+}} ( \yzero ) = 0$.
It now follows from Theorem~1 that $\vv\not\in Y_{0}$, and thus
that $h= \hzero$ and $\yzero$ is optimal.~$\quad\Box$
\begin{cor}
The components of a best concave approximation to data $f_{j}$,
$r \leq j \leq s$, are given from Algorithm 1 with $\sigma=-$ by
$y_{j}=f_{j}(I)-h$, where in this case $I=I^{-}(r,s)$ and
$h=h^{-}(r,s)$.
\end{cor}

\subparagraph{}

The algorithms in the next subsections will join optimal vertex
sets produced by Algorithm 1 of consecutive ranges of the data,
and they will also cut such optimal vertex sets in two. It is
convenient to prove here that the resulting sets remain optimal
vertex sets. The proof requires one further property of the
optimal vertex sets produced by Ubhaya's algorithm.

This algorithm is a specialization of an algorithm by Graham
(1972) for finding the convex hull of an unordered set of points
in the plane. It is logically equivalent to the algorithm of
Kruskal (1964) for monotonic approximation, which is more
efficient than the algorithm of Miles (1959) for monotonic
approximation, but produces the same results. Two conclusions
follow from this. The first, which will be needed later, is that
if $i$ and $k$ are consecutive elements of \isig{r}{s}\ then
    \beqa
 \sigma   g_{ik} ({\bf f}(I)) &=& \min
        \{ \, \sigma g_{ij}({\bf f}(I)) \, : \, i \leq j    \leq s  \,
        \} \qquad \mbox{ and } \label{2110a} \\
  \sigma  g_{ik}({\bf f}(I))&=& \min
        \{ \, \sigma g_{jk} ({\bf f}(I))\, : \, r \leq j    \leq i  \,
        \}.
    \label{2110b}
    \eeqa

\subparagraph{}

The second, which is of some theoretical interest, is that the
gradients $g_{j \,j+1}({\bf f}(I))$, for $1 \leq j \leq n-1$, are
the best \emph{least squares} approximations to the numbers
 $g_{j \,j+1}(\ff)$, for $1 \leq j \leq n-1$, subject to to the
constraint that they monotonically increase, {\em i.e.} to
convexity of the points of which they are gradients. The proof of
this interesting equivalence is a straightforward consequence of
the equivalence of Kruskal's and Miles's algorithms for monotonic
approximation. It is given in Cullinan (1986).

The results for the joining and splitting of optimal vertex sets
require definition of trivial optimal vertex sets by
\beq
    \isig{r}{s} = [r,s] \qquad \mbox{  when } s-1 \leq r \leq s .
    \label{2113}
\eeq
It is also convenient to define $\hsig{r}{s} = 0$ when $r \geq s-1$.
The following lemmas then hold.

\begin{lemma}
A subset of one or more consecutive elements of an optimal vertex
set is itself an optimal vertex set.
\end{lemma}
\paragraph{Proof} The proof is trivial unless the subset has at
least four elements, when it follows from the nature of extreme
points of convex sets using (\ref{2110a}) and (\ref{2110b}).~$\quad
\Box$

\subparagraph{}

The second lemma gives conditions for the amalgamation of optimal
vertex sets.
\begin{lemma}
Given vertex sets \isig{r}{s}\ and \isig{s}{t}\, necessary and
sufficient conditions for
    \beq
       \isig{r}{t} = \isig{r}{s} \union \isig{s}{t}
    \label{2111}
    \eeq
are that there exist $r' \leq r$, $r' \in \isig{r}{s}$, and $t'
\geq t$, $t' \in \isig{s}{t}$, such that
    \beq
    s \in \isig{r'}{t'}.
    \label{2112}
    \eeq
\end{lemma}
\paragraph{Proof} The proof is trivial except for the
sufficiency of (\ref{2112}) when $r<s<t$. Let $s^{-}$ be the left
neighbour of $s$ in $\isig{r}{s}$ and $s^{+}$ its right neighbour
in $\isig{s}{t}$.
It follows from (\ref{2110a}) and (\ref{2110b}) that (\ref{2111}) holds provided that
    \[
  \sigma  c_{s^{-},s,s^{+}}(\ff ) > 0.
    \]
However, it follows from (\ref{2112}) that
    $\sigma c_{isj}(\ff ) > 0$
for \emph{all} $i$ and $j$ in the range $r' \leq i<s<j \leq t'$,
and this range contains $s^{-}$ and $s^{+}$.~$\quad \Box$
%
%
\subsection{The case $q=1$}
The algorithm to be presented constructs ranges $[1,s]$ and
$[t,n]$, where $s \leq t$, a price $h \geq 0$, and a vertex set
$I$ of $[1,n]$ such that
\[
I = I^{+} (1,s) \union I^{-} (t,n) .
\]

The best approximation \yone\ is then given by the final value of
the vector \yy\ defined by
\beqa
    y_{i} & = & f_{i} + h \qquad \mbox{ when } i \in I^{+}(1,s)
      \label{2201} \\
    y_{i} & = & f_{i} - h \qquad \mbox{ when } i \in I^{-}(t,n)
        \label{2202} \\
    y_{j} & = & y_{j} (I) \qquad 1 \leq j \leq n
    \label{2203}
\eeqa

It will be shown that $s=t$ only if $h=0$, so that \yone\ is well
defined. This construction is illustrated in Figure~2. Note that when $s=t=n$, $\yy = \yzero$.
%
%
\begin{figure}
\begin{center}
\begin{picture}(504,216)
\setlength{\unitlength}{0.75pt}
\thicklines
\dosmblob{36}{180}
\writeat{36}{189}{$y^{1}_{1}$}
\doline{36}{180}{1,-2}{36}
\dosmblob{72}{108}
\doline{72}{108}{3,-2}{54}
\dosmblob{126}{72}
\doline{126}{72}{4,1}{72}
\dosmblob{198}{90}
\writeat{198}{99}{$y^{1}_{s}$}
\dosmblob{306}{126} \writeat{306}{135}{$y^{1}_{t}$}
\doline{306}{126}{4,1}{72} \dosmblob{378}{144}
\doline{378}{144}{1,0}{54} \dosmblob{432}{144}
\doline{432}{144}{1,-2}{36} \dosmblob{468}{72}
\writeat{468}{92}{$y^{1}_{n}$} \thinlines
\doline{198}{90}{3,1}{108} \writeat{234}{126}{join}
\doblob{36}{140} \writeat{36}{128}{1} \doline{36}{140}{1,-2}{36}
\doblob{72}{68} \doline{72}{68}{3,-2}{54} \doblob{126}{32}
\doline{126}{32}{4,1}{72} \doblob{198}{50} \writeat{198}{38}{$s$}
\doblob{306}{166}
\writeat{306}{178}{$t$}
\doline{306}{166}{4,1}{72}
\doblob{378}{184}
\doline{378}{184}{1,0}{54}
\doblob{432}{184}
\doline{432}{184}{1,-2}{36}
\doblob{468}{112}
\writeat{468}{127}{$n$}
\end{picture}
\end{center}
\caption{Construction of best convex-concave approximation}
\end{figure}
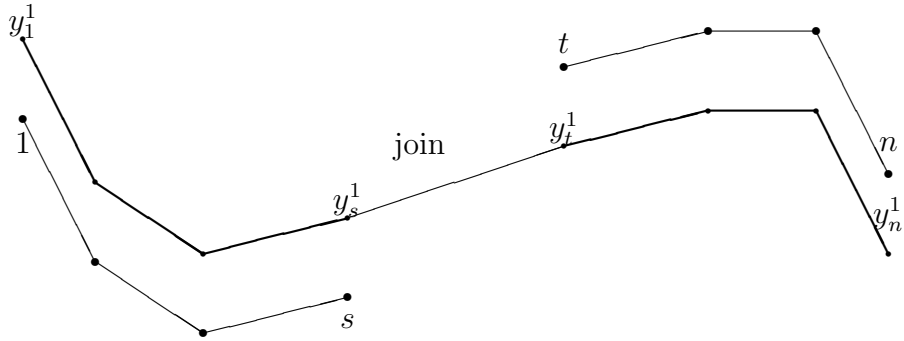
%
%

The algorithm to be given is believed to be more efficient than
that in Cullinan \& Powell (1982). For example if the data are in
$Y_{0}$, the new algorithm will require only iteration,
whereas the former one only has this property if the data lie on
a straight line.

The algorithm builds up $I$ by looking alternately at the left
and right ranges. Beginning with $I=\{1,n \}$, $s=1$, and $t=n$,
it adds an index $k$ of \iplus{s}{t}\ to $I$ if the least possible
final value of $F$ consistent with doing this is not greater than
the least possible value of $F$ consistent with ending the
calculation with the existing value of $s$. After adding one
index in \iplus{s}{t}\ to $I$ and increasing $s$ to the value of this index, it then examines the next. When it
is not worth adding any more indices from \iplus{s}{t}\ to $I$,
it tries to add indices in \iminus{s}{t}\ to $I$ working
backwards from $t$, adding $k$ to $I$ if it is not necessarily
more expensive to finish with $t$ reduced to $k$ than with $t$ as
it is and decreasing $t$. After indices have been added to $I$ from \iminus{s}{t}\
it may then be possible to add more to $I$ from the new $\iplus{s}{t}$, so
the process alternates between \iplus{s}{t}\ and \iminus{s}{t}\
until $s$ equals $t$ or $t-1$ or until the algorithm fails twice
running to add any indices.

\begin{alg}
To find a best convex-concave approximation when $n \geq 2$.
\end{alg}

\balg
\item Set $s:=1$, $t :=n$, $h := 0$, and $I := \{ \, 1,n \, \}$.

\item If $s \geq t-1$: {\bf stop.} Otherwise: set $u = t-s$ and $I' := \iplus{s}{t}$.

\item Let $s^{+}= s^{+} (I')$. Find $h^{+}(s,s^{+})$ and set
$h' := \max( h, h^{+} (s,s^{+}))$. \\
If $f_{t}(s,s^{+}) > f_{t} - 2 h'$: go to Step~5.

\item Add $s^{+}$ to $I$ and delete $s$ from $I'$.
Set $h:= h'$ and $s := s^{+}$. \\
If $s < t-1$: go back to Step~3.

\item If $s = t$: {\bf stop.} Otherwise: set $I' := I^{-}(s,t)$.

\item Let $t^{-} = t^{-} (I')$. Find $h^{-}(t^{-},t)$ and set $h' := \max (h, h^{-}(t^{-},t))$. \\
If $f_{s}(t^{-},t) < f_{s} + 2 h'$: go to Step~8.

\item Add $t^{-}$ to $I$ and delete $t$ from $I'$.
Set $h:= h'$ and $t := t^{-}$. \\
If $s < t $: go back to Step~6.

\item If $t-s = u$: {\bf stop.} Otherwise: go back to Step 2. \ealg

The optimal vertex sets $I'$ are found using Algorithm~1 or
Corollary 2.1. The prices \hsig{s}{t}\ are found from
(\ref{2101}) and $f_{t}(s,s^{+})$ and $f_{s}(t^{-},t)$ from
(\ref{2103b}).

The first remark is that $s$ is non-decreasing, $t$ is
non-increasing, and $s \leq t$ with $s=t$ only when $h=0$. The
vector \yone\ is therefore well-defined by
(\ref{2201})--(\ref{2203}).

The possibility that Algorithm~1 can end with $s=t$ will involve slightly more complexity
when this Algorithm is used later on when $q \geq 3$.
It might be prevented by relaxing either of the inequalities in Steps~3and~6.
However if {\em both} these inequalities are relaxed then the Algorithm will fail.
For example with equally spaced data and $\ff = ( \, 0,1,0,1 )$, relaxing both inequalities yields
$s=1$, $t=4$, and $F( \yy ) = \frac{2}{3}$, whereas the Algorithm as it stands correctly calculates
$h = \half$ and $\yy = ( \, \half , \half , \half , 1  )$.
It seemed best to let both inequalities be strict for reasons of symmetry.

The next result concerns the conditions under which the Algorithm
increases $s$ and decreases $t$.

\begin{lemma}
At any entry to Step~2 of Algorithm~2, let the vector \yy\ be
defined by (\ref{2201})--(\ref{2203}). Then the Algorithm
strictly decreases $t-s$ if and only if \beq
    \mbox{there exists } j \: : s<j<t \quad \mbox{ and } \quad
    |f_{j} - y_{j}| \geq h .
\label{151}
\eeq
\end{lemma}
(This somewhat cumbersome statement is needed to include the case
where $h=0$.)

Since the points $(x_{j},y_{j})$, $s \leq j \leq t$, are
collinear, (\ref{151}) is equivalent to the statement that there
exists a point of the graph of the data between $x_{s}$ and
$x_{t}$ lying on or outside the parallelogram $\Pi(h)$ with
vertices $(x_{s},f_{s})$, $(x_{s},f_{s}+2h)$, $(x_{t},f_{t})$, and
$(x_{t},f_{t}-2h)$. This parallelogram $\Pi(h)$ is illustrated in
Figure~3.
%
%
\begin{figure}
\begin{center}
\begin{picture}(180,180)
\thicklines
\doblob{36}{36}
\writeat{36}{27}{$s$}
\doline{36}{36}{2,1}{144}
\doline{180}{108}{0,1}{36}
\writeat{188}{117}{$h$}
\writeat{188}{136}{$h$}
\doblob{180}{144}
\writeat{180}{153}{$t$}
\doline{36}{36}{0,1}{36}
\writeat{28}{44}{$h$}
\writeat{28}{64}{$h$}
\doline{36}{72}{2,1}{144}
\writeat{36}{80}{$\Pi(h)$}
\thinlines
\dosmblob{36}{54}
\doline{36}{54}{2,1}{144}
\dosmblob{180}{126}
\end{picture}
\end{center}
\caption{The parallelogram $\Pi(h)$}
\end{figure}
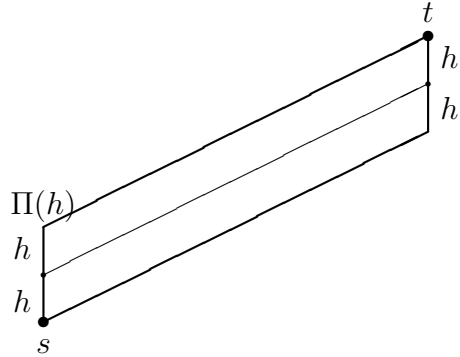
%
\paragraph{Proof}

In the trivial case when $s \geq t-1$, (\ref{151}) is false and
the Algorithm stops without altering $s$ or $t$.

Otherwise, there are one or more data points $f_{j}$ with $s<j<t$.
Suppose first that (\ref{151}) does not hold. Then $h>0$. (If
$h=0$ then (\ref{151}) holds trivially!) It will be shown that in
this case both $s$ and $t$ are left unchanged. Step~3 will
calculate an index $s^{+} : s<s^{+} \leq t$.

If $s^{+} = t$, then because $h' \geq h>0$, it immediately
follows that $f_{t}(s,s^{+}) = f_{t} > f_{t} - 2h'$, so that
Step~3 will lead immediately to Step~5 and $s$ will not be
increased.

If $s^{+} < t$, let $y_{s^{+}}$ be defined by
(\ref{2201})--(\ref{2203}) and let $g= g_{st}(\yy)$. Then by
hypothesis, $f_{s^{+}}> y_{s^{+}} - h$, which implies that $g_{s
\; s^{+}} (\ff ) > g$. Thus,
\begin{eqnarray*}
f_{t} (s,s^{+}) & = &  f_{s} +  g_{s \, s^{+}} (\ff )(x_{t}-x_{s}) \\
           & > &  f_{s} + g ( x_{t} - x_{s}) = y_{t}-h = f_{t}-2h
           \\
           & \geq &  f_{t} - 2h' .
\end{eqnarray*}
Therefore Step~3 will again lead immediately to Step~5.

The same arguments applied to $t^{-}$ calculated by Step~6 show
that Step~6 will branch immediately to Step 8 and so $t$ will also
be left unchanged.

\subparagraph{}

Now suppose conversely that (\ref{151}) does hold, so that there
is an index $j$ with data point lying on or outside $\Pi$.
It will be shown that in this case if the algorithm does not increase $s$, it must then
decrease $t$.

Consider first the case where $f_{j} \leq y_{j} -h$.
Step~3 will calculate an index $s^{+}$ in the range $s < s^{+}
\leq t$, and, from (\ref{2110a}), with $g_{s \; s^{+}} (\ff )
\leq g_{sj} (\ff )$. Then \beqa
 f_{t}(s,s^{+})  & = & f_{s} + g_{s \; s^{+}} (\ff ) (x_{t}-x_{s}) \\
            & \leq & f_{s} + g_{sj}(\ff )(x_{t}-x_{s}) \\
            & = & f_{s} + ((f_{j} - f_{s})/(x_{j} - x_{s}))(x_{t}-x_{s}) \\
            & \leq & f_{s} + ((y_{j}-h - f_{s})/(x_{j} - x_{s})) (x_{t}-x_{s})
                    \mbox{ by hypothesis} \\
            & = & y_{s} - h + ((y_{j}-y_{s})/(x_{j}-x_{s}))(x_{t}-x_{s})\\
            & = & y_{t} - h, \quad \mbox{ by the definition of }y_{j}.
\eeqa
It now follows immediately that when $h' = h$, Step~4 will be entered and
$s$ increased to $s^{+}$, as required.

Now suppose that $h' = h^{+}(s,s^{+}) > h$ and that $s$ is not increased, i.e. that the test in
Step~3 leads to Step~5. Then $f_{t}(s,s^{+}) > f_{t} - 2h' =
f_{t} - 2h^{+}(s,s^{+})$. Let $h^{+}=h^{+}(s,s^{+})$. By
definition of $h^{+}$, there must be an index $i$ such that
$f_{i} - f_{i}(s,s^{+}) = 2h^{+}$.
This case is illustrated in Fig.~4.
It demonstrates the heart of the principle behind the algorithm, because the four
data points with indices $s$, $i$, $s^{+}$, and $t$ determine a
lower bound on $\inf F(Y_{1})$.
%
%
\begin{figure}
\begin{center}
\begin{picture}(180,180)
\thinlines
\doblob{36}{36}
\writeat{36}{27}{$s$}
\doline{36}{36}{2,1}{144}
\doline{180}{108}{0,1}{36}
\doblob{180}{144}
\writeat{180}{153}{$t$}
\doline{36}{36}{0,1}{36}
\doline{36}{72}{2,1}{144}
\writeat{24}{72}{$\Pi(h)$}
\thicklines
\doline{36}{36}{6,1}{144}
\doblob{108}{48}
\writeat{108}{39}{$s^{+}$}
\doline{180}{60}{0,1}{84}
\writeat{192}{102}{$2d$}
\doline{36}{36}{0,1}{84}
\thinlines
\doline{36}{120}{0,1}{48} \doblob{72}{144}
\writeat{81}{144}{$i$} \doline{72}{42}{0,1}{102}
\writeat{60}{93}{$2 h^{+} $} \doline{36}{144}{1,0}{144}
\doblob{132}{152}
\writeat{136}{161}{$t^{-}$}
\doline{180}{144}{-6,1}{144}
\doblob{156}{90}
\writeat{164}{90}{$k$}
\doline{156}{148}{0,-1}{58}
\writeat{166}{119}{$2 h'$}
\end{picture}
\end{center}
\caption{When $s$ does not increase, $t$ must decrease}
\end{figure}
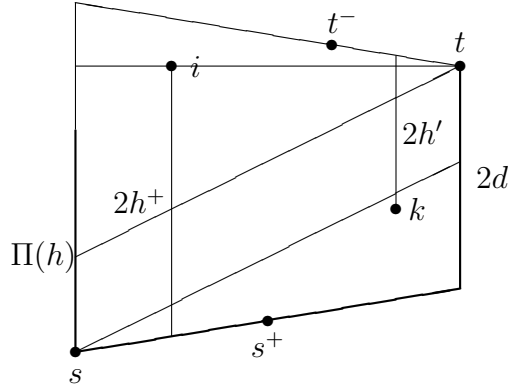
%
%

Define \beq
    2d = f_{t}- f_{t}(s,s^{+}).
    \label{defofd}
\eeq
Then by the hypothesis that Step~3 led to Step~5, \beq
    h^{+} > d . \label{beta2}
\eeq Step~5 will be entered with $h<h^{+}$ and will calculate
$I^{-}(s,t)$. Step~6 will then find an index $t^{-} : s \leq
t^{-} <t$. It follows from (\ref{2110a}) that $f_{i} \leq
f_{i}(t^{-},t)$, and so \beq 2h^{+}= f_{i} - f_{i} (s,s^{+}) \leq
f_{i}(t^{-},t) - f_{i} (s,s^{+}). \label{1700} \eeq This
inequality along with (\ref{defofd}) and (\ref{beta2}) then show
that the function $l \mapsto f_{l} (t^{-},t) - f_{l} (s,s^{+})$
is strictly decreasing.

If $h$ does not increase in Step~7, it follows immediately, from
the assumption that $h<h^{+}$, that $f_{s} (t^{-},t) - f_{s}
(s,s^{+})= f_{s}(t^{-},t) - f_{s} > 2h^{+} > 2h$, so that $t$
will be reduced to $t^{-}$ by the test in Step~6.

If, on the other hand, $h$ does increase in Step~7, its new value
$h'$ is given by $2h' = f_{k}(t^{-},t) - f_{k}$ for some $k$ in
the range $s<k<t$. But the definition of $s^{+}$ implies that
$f_{k} \geq f_{k}(s,s^{+})$, and so $2h' \leq f_{k}(t^{-},t) -
f_{k}(s,s^{+}) < f_{s}(t^{-},t) - f_{s}(s,s^{+})$, so that again
$t$ must be reduced to $t^{-}$ by Step~7.
In fact it is easy to show that in this case, $t^{-} \geq i$, $h' < h^{+}$, and that
$\inf F(Y_{1}) \geq \half (f_{s^{+}}(i,t)- f_{s^{+}})$.

Thus $t$ will be reduced to $t^{-}$ in all cases.

In the case where there exists $j$ such that $f_{j} \geq y_{j} +
h$, if $s$ is not increased immediately after the next entry to
Step 3, the same argument shows that either $t$ must be reduced
in the next operation of Steps~6 and 7, or $s$ must be increased
immediately thereafter. ~$\quad \Box$
%
%
\begin{lemma}
At any exit from Step~4 or   Step~7, let \yy\ be defined by
(\ref{2201})--(\ref{2203}). Then \beq \mbox{when $s>1$,} \quad
\civ{s-1}{y} \geq 0, \label{23012} \eeq and \beq \mbox{when $t<n$,
} \quad \civ{t-1}{y} \leq 0. \label{23013} \eeq
\end{lemma}
\subparagraph{Proof} Suppose first that Step~3 is entered. Then
(\ref{23012}) is equivalent to \beq c_{s^{+}-1}( \yy ) \geq 0 ,
\eeq which is in turn equivalent to the identity \beq f_{t} -
f_{t}(s,s^{+}) \geq 2 h' . \label{2343}\eeq But this is simply
the test leading to Step~4.

If $t<n$, there will exist an index $t^{+}=t^{+}(I)$, and
(\ref{23013}) will be equivalent to  the inequality \beq
f_{s^{+}}(t,t^{+}) - f_{s^{+}} \geq 2 h' . \label{2341} \eeq

Suppose firstly that $h'>h$, so that there exists $i$ such that
$s<i<s^{+}$ and \beq f_{i} - f_{i}(s,s^{+}) = 2 h' . \label{2342}
\eeq

The earlier definition of $t$ from $t^{+}$ implies that $f_{i} \leq
f_{i}(t,t^{+})$
and it follows from (\ref{2342}) and (\ref{2343})
that the monotonic function $\phi : l \mapsto f_{l}(t,t^{+}) -
f_{l}(s,s^{+})$ satisfies $\phi (i) \geq 2h'$ and $\phi(t) \geq 2
h'$. It follows that $\phi(s^{+}) \geq 2 h'$, which is equivalent
to (\ref{2341}).
It follows that whenever $h$ increases, both the join constraints
$c_{s-1}$ and $c_{t-1}$ (when defined) are feasible.

When Step~4 is entered and $h$ does not increase, if $c_{t-1}$ is
initially non-positive, which is equivalent to the inequality \beq
f_{s}(t,t^{+}) - f_{s} \geq 2h , \eeq it follows from this and
(\ref{2343}) with $h'=h$ that $\phi (s) \geq 2h$ and $\phi (t)
\geq 2h$ so that $\phi (s^{+}) \geq 2 h$ as required.

Similarly when Step~7 is entered and $h$ does not increase, if
$c_{s-1}$ is initially non-negative, it will remain so when $s$
is increased.

The result then follows by induction. ~$\quad \Box$

\begin{cor}
If steps~2 to 8 of Algorithm~2 are executed with $h$ initially set to any number $\hb > 0$, then
(\ref{23012}) and (\ref{23013}) remain satisfied.
\end{cor}

The effectiveness of Algorithm~2 will now be established.
\begin{theorem}
Algorithm~2 produces integers $s$ and $t$ with $s \leq t$, an
index set $I$ such that \beq I = I^{+} (1,s) \union   I^{-} (t,n)
, \label{23one} \eeq and a real number $h$ such that $h=0$ if and
only if $s=t$. If \yone\ is then defined by
(\ref{2201})--(\ref{2203}), then $F( \yone ) = h = \inf F(Y_{1})$
and $\yone \in Y_{1}$.
\end{theorem}
\paragraph{Proof}
The first stage is to establish (\ref{23one}). Assume inductively
that it holds before a series of sections is added, and without
loss of generality that a series of convex sections with indices
\stoalpha\ are added by successive entries to Step~4 for the same
value of $t$. Each time an index is deleted from $I'$ it follows
immediately from Lemma~1 that the new value of $I'$ is also an
optimal index set, so it always holds that $I' = I^{+}(s,t)$. It
also follows from Lemma~1 that $\stoalpha = I^{+}(s,s_{\alpha})$.
Let $I_{1} = \{ i \in I \: : i \leq s \}$. After all the sections
are added, $I_{1} = I^{+}(1,s) \cup I^{+}(s,s_{\alpha})$. If
$s=1$, it follows at once from Lemma~2 that
$I_{1}=I^{+}(1,s_{\alpha})$ as required. Otherwise, let $t'$ be
the value that $t$ had when $s$ was increased from $s^{-}$. Then
$s_{\alpha} \leq t \leq t'$. Immediately after $s$ was increased from $s^{-}$,
$s \in I^{+} (s^{-},t')$. The conditions of Lemma~2 are
therefore satisfied and $I_{1} = I^{+}(1,s_{\alpha})$. The same
argument applied to concave sections then establishes (\ref{23one}).

The algorithm clearly produces a number $h$ such that $s=t$ only
if $h=0$. It is a consequence of Lemma~3 that if $s< t-1$ and
$h=0$, then the algorithm will reduce $t-s$.  If $s=t-1$ and
$h=0$, Step~4 will increase $s$ to $t$. Thus $h=0$ if and only if
$s=t$.  Thus \yone\ is well defined by (\ref{2201})--(\ref{2203}).
The number $h$ is given by $h= \max (h_{(1)},h_{(2)})$ where \beq
h_{(1)} = h^{+} (1,s) , \label{23seven} \eeq \beq h_{(2)} =
h^{-}(t,n) . \label{23eight} \eeq It follows from (\ref{23one})
and Lemma~3 that when the algorithm terminates, \beq F(\yone )=h.
\label{23021} \eeq

It follows directly from (\ref{23one}), (\ref{2201})--(\ref{2203})
and Lemma~4 that $\yone \in Y_{1}$.

It remains to prove that \yone\ is optimal. The method of proof
chosen to do this can be simplified in this case, but generalizes
more directly to the case of $q > 1$. If $h=0$, optimality follows
immediately from (\ref{23021}). Otherwise suppose that there is a
vector \vv\ such that $F( \vv ) < h$.
The price $h$ will be defined from
(\ref{2101})
with particular values of $\sigma$, $r$, and $s$.
Let $\jst$ \label{jstardef}
be the lowest value of $j$ in
this equation that defines the final value of $h$.
Then \jst\ lies strictly
between two neighbouring elements $k$ and $k^{+}$ of $I$. Assume
firstly that $\jst <s$. Define the set $K = \{ \, k, \jst, k^{+},s,t \,
\}$. Since $h>0$, $s<t$, so $K$ has at least four elements (it is
possible that $k^{+}=s$). Now \bea
y^{1}_{k} & = & f_{k} + h,  \\
y^{1}_{ \jst } & = & f_{ \jst } - h,  \\
y^{1}_{k^{+}} & = & f_{k^{+}} + h,  \\
y^{1}_{s} & = & f_{s} + h,  \\
y^{1}_{t} & = & f_{t} - h. \eea
Since $F( \vv ) < h$, \bea
v_{k} & < & f_{k} + h,  \\
v_{ \jst } & > & f_{ \jst } - h,  \\
v_{k^{+}} & < & f_{k^{+}} + h,  \\
v_{s} & < & f_{s} + h,  \\
v_{t} & > & f_{t} - h. \eea
By definition of $y^{1}_{ \jst }$,
\[
\civ{k \jst k^{+}}{\yone } = 0,
\]
and so $\civ{k \jst k^{+}}{\vv} < \civ{k \jst k^{+}}{\yone } = 0$.
From Lemma~4, $\civ{i-1}{\yone } \geq 0$ for all $i$ in the range
$ \jst \leq i \leq t-1$. It follows from the corollary to Theorem~1 that
\[
\civ{ \jst s t}{\yone } \geq 0.
\]
Then $\civ{ \jst s t}{\vv} > \civ{ \jst s t}{\yone} \geq 0$. If
$k^{+}=s$, Theorem~1 can be immediately applied to $K$ to  show
that $\vv \not\in Y_{1}$. If $k^{+}<s$, it follows from the corollary to Theorem~1
that because $\civ{ \jst s t}{\vv}>0$, at least one of the
consecutive divided differences $\civ{ \jst k^{+} s}{\vv}$ and
$\civ{k^{+} s t}{\vv}$ must be positive, so that again $\vv
\not\in Y_{1}$.

If $ \jst >s$, let $K= \{ \, s,t,k, \jst , k^{+} \, \}$. Then the same
argument shows that \vv\ cannot be feasible. Therefore \yone\ is
optimal.~$\quad \Box$
%
%
%
\subsection{The Case $q=2$}
The best $Y_{1}$ approximation constructed in Section~2.2 was
defined by equations (\ref{2201})--(\ref{2203}) from the pieces
$[1,s]$ and $[t,n]$, the index set $I$, and the price $h$. The
best $Y_{2}$ approximation will in general be constructed from
three pieces $P_{1}=[1,s_{1}]$, $P_{2}=[t_{1},s_{2}]$,
$P_{3}=[t_{2},n]$, a price $h \geq 0$ with $h>0$ only when
$s_{1}<t_{1}$ and $s_{2}<t_{2}$, and an index set $I$ of $[1,n]$
such that \beq I= I^{+}(P_{1}) \union I^{-}(P_{2}) \union
I^{+}(P_{3}), \label{240300} \eeq as the ultimate value ${\bf
y}^{2}$ of the vector \yy\ defined by the equations \beqa y_{i} &
= & f_{i} +  (-)^{ \alpha - 1}h \qquad \mbox{ when } i \in
I \cap P_{ \alpha }, \qquad  1 \leq \alpha \leq 3,\label{240301} \\
 y_{j} & = & y_{j} (I) \qquad 1 \leq j \leq n. \label{240303}
\eeqa

%
%
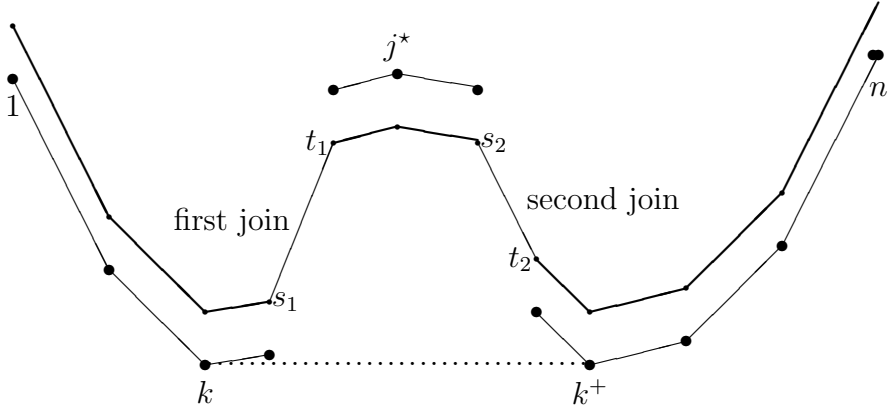
\begin{figure}
\begin{center}
\begin{picture}(360,156)
\thinlines
\doblob{18}{134}
\writeat{18}{124}{$1$}
\doline{18}{134}{1,-2}{36}
\doblob{54}{62}
\doline{54}{62}{1,-1}{36}
\doblob{90}{26}
\writeat{90}{16}{$k$}
\doline{90}{26}{6,1}{24}
\doblob{114}{30}
\multiput(90,26)(4,0){36}{.}
\doblob{234}{26}
\writeat{234}{16}{$k^{+}$}
\doline{234}{26}{-1,1}{20}
\doblob{214}{46}
\doline{234}{26}{4,1}{36}
\doblob{270}{35}
\doline{270}{35}{1,1}{36}
\doblob{306}{71}
\doline{306}{71}{1,2}{36}
\doblob{342}{143}
\writeat{342}{130}{$n$}
\doblob{162}{136}
\writeat{162}{146}{$\jst$}
\doline{162}{136}{6,-1}{30}
\doblob{192}{130}
\doline{162}{136}{-4,-1}{24}
\doblob{138}{130}
\thicklines
\dosmblob{18}{154}
\doline{18}{154}{1,-2}{36}
\dosmblob{54}{82}
\doline{54}{82}{1,-1}{36}
\dosmblob{90}{46}
\doline{90}{46}{6,1}{24}
\dosmblob{114}{50}
\writeat{120}{50}{$s_{1}$}
\dosmblob{234}{46}
\doline{234}{46}{-1,1}{20}
\dosmblob{214}{66}
\writeat{208}{66}{$t_{2}$}
\doline{234}{46}{4,1}{36}
\dosmblob{270}{55}
\doline{270}{55}{1,1}{36}
\dosmblob{306}{91}
\doline{306}{91}{1,2}{36}
\doblob{340}{143}
\dosmblob{162}{116}
\doline{162}{116}{6,-1}{30}
\dosmblob{192}{110}
\writeat{198}{110}{$s_{2}$}
\doline{162}{116}{-4,-1}{24}
\dosmblob{138}{110}
\writeat{132}{110}{$t_{1}$}
\thinlines
\doline{114}{50}{2,5}{24}
\writeat{100}{80}{first join}
\doline{192}{110}{1,-2}{22}
\writeat{239}{88}{second join}
\end{picture}
\end{center}
\caption{Construction of $Y_{2}$ approximation}
\end{figure}
%
%
The construction of \yy\ is illustrated in Figure~5.

The Algorithm will construct this best $Y_{2}$ approximation from
the quantities $h_{0}=h^{+}(1,n)$ and $I^{+}(1,n)$ constructed by
Algorithm~1. If the value $h_{0}$ of this approximation is zero,
the best approximation over $Y_{0}$ is also a best approximation
over $Y_{2}$. Otherwise, $h_{0}> 0$ is determined by three data
$f_{k}$, \fjs , and $f_{k^{+}}$ such that $k$ and $k^{+}$ are
consecutive elements of $I^{+}(1,n)$ and $k< \jst < \kpl$. The
discussion in Section~2.1 shows that unless the divided
differences of \yy\ change sign at least once in the range $[ \,
x_{k} , \, \sbkp{x} \, ]$, then $F( \yy ) \geq h_{0}$. The set
$I^{-}(t_{1},s_{2})$ is therefore put in this range. The
Algorithm begins with $s_{1}=k$, $t_{1}=s_{2}=\jst$, and
$t_{2}=k^{+}$. It then sets $h= \max
(h^{+}(1,s_{1}),h^{+}(t_{2},n))$. Next, it uses Algorithm~2
modified to calculate a best convex-concave approximation to the
data with indices in $[1, \jst ]$ consistent with paying this
minimum price $h$. It is an important feature of the problem that
this can be done by starting Algorithm~2 with $s=k$ and $t=
\jst$, i.e. the existing elements of $I^{+}(1, \jst)$ below
$k$ can be kept in place.

This process can increase $s_{1}$ beyond $k$, reduce $t_{1}$
below \jst , and increase $h$. Let $h^{(1)}$ be its new value. A
best {\em concave-convex} approximation to $\fjs , \ldots , f_{n}$
starting with $h=h^{(1)}$ is next identified by applying a
modified version of Algorithm~\thealg\ with $s= \jst$ and $t=
\kpl$, in general increasing $s_{2}$ beyond \jst\ and reducing
$t_{2}$ below $k^{+}$. If this second calculation does not
increase $h$ above $h^{(1)}$, the best approximation \yy\ can be
constructed immediately from (\ref{240301})--(\ref{240303}). If,
however, $h>h^{(1)}$, there is the complication that the lower
value of $h$ when the first calculation took place may have
joined too many sections for the first join constraints $c_{s_{1}
-1}(\yy )$ and $c_{t_{1}-1}(\yy )$ determined by the new higher
value of $h$ to have the right signs. In this case the remedy
proposed is to repeat the first calculation starting with the new
value of $h$.

The following algorithms will therefore  require a modified
version of Algorithm~\thealg\ to carry out a best convex-concave
approximation or a best concave-convex approximation on a range
of data, starting with a prescribed value of $h$. This task is
best carried out by modifying Algorithm~\thealg\ in detail, but
the following description is equivalent and simpler. The
following procedure calculates a best approximation to the data
in the range $[ \sal , \tal ]$ compatible with an existing price
$h$, the approximation being over $Y_{1}$ or $Y_{-1}$ according
as $\alpha$ is odd or even. It can be seen as trying to close the
join $[ \sal , \tal ]$ as much as possible by constructing the
best approximation to this range of data compatible with the
given starting value of $h$.

The notation for the pieces is is based on the observation that
when $h>0$, $\tal = s^{+}_{\alpha}(I)$ and when $h=0$, $\tal =
\sal$, so that the location of each of the two joins
$[s_{1},t_{1}]$ and $[s_{2},t_{2}]$ can be specified simply by
consecutive elements of $I$. Thus the location of the pieces and
joins can be specified simply through the quantities $s_{\alpha
}$, $1 \leq \alpha \leq 2$. It is convenient to regard these as
members of an ordered subset $S$ of $I$, and to include $n$ in
$S$. Given an index set $I$, define a {\em piece set } $S$ of $I$
to be an ordered subset of $I$ such that $n \in S$.
Then $\sal$ will denote the $\alpha$th element of $S$.

\begin{alg}
closejoin($\alpha$)
\end{alg}
\balg
\item If $\alpha$ is even: replace \ff\ by $- \ff$. \\
      Set $s = \sal$ and $t = s^{+}_{ \alpha }(I)$.
\item Carry out Steps~2 to~8 of Algorithm~2.
\item Set $s_{ \alpha } = s$. \\
If $\alpha$ is even: replace \ff\ by $- \ff$. \ealg Note that
this procedure modifies $h$, $I$, and $S$.

The algorithm constructs $\yy ^{2}$ by calculating the appropriate
pieces and a price $h(S)$ from which $\yy ^{2}$ is constructed.
The notation used for keeping track of the pieces needs to cover
the case of pieces that consist of only one point, and to
generalize easily when $q>2$. It also has to cope with the
trivial cases where $\ff \in Y_{0}$ or $\ff \in Y_{1}$ and there
are, therefore, only one or two pieces instead of three.

Given a piece set $S$, define its price $h(S)$ as follows. Let
$q'=|S|-1$, $t_{0}=1$, and when $q' \geq 1$ define $t_{ \alpha }
= s^{+}_{ \alpha } (I)$ when $ 1 \leq \alpha \leq q'$. Then
\beq h(S) = \max \{ \; h^{(-)^{ \alpha -1}}( \talm, \sal
) \: : \: 1 \leq \alpha \leq q'+1 \; \}. \label{240305} \eeq
Recall that $h^{\sigma}(t,s) = 0$ whenever $s \geq t-1$.

It is also worth recording the location of the data points
determining the optimal value of $h$. Given a piece set $S$, when
$h(S)>0$ it follows from (\ref{240305}) and (\ref{2101}) that
there exists a lowest index $\jst(S)$ such that \beq h(S) = \half
\sigma ( f_{ \jst } - f_{ \jst } (k, k^{+})),
\label{10601} \eeq where $\sigma = (-)^{ \beta -1} $, $1 \leq \beta
\leq q'+1$, and \beq t_{ \beta -1} \leq k < \jst < k^{+} \leq s_{
\beta} . \label{20601} \eeq Once $\jst$ is known, the quantities
$k$, $k^{+}$,and $\beta$ are uniquely determined by (\ref{10601})
and (\ref{20601}).

The following algorithm then constructs $h$, $I$ and $S$
determining $\yy ^{2}$.

\begin{alg}
To find a best approximation over $Y_{2}$.
\end{alg}
\balg
\item   Set $S= \{n\}$. Use Algorithm~1 to calculate $I=I^{+}(1,n)$ and $h=h^{+}(1,n)$. \\
    If $h=0$: {\bf stop.} \\
    Otherwise: find $\jst, k, k^{+}$ determined from $S$ by (\ref{10601}) and proceed to Step~2.
\item   Insert $\jst$ into $I$ and $k$, $\jst$ into $S$. Calculate $h = h(S)$.
\item Apply {\em closejoin}(1). Set $h^{(1)}=h$.
\item Apply {\em closejoin}(2). If $h=h^{(1)}$: {\bf stop.} \\
    Otherwise proceed to Step~5.
\item If $h^{(1)}=0$: go to Step~6. \\
    If $s_{1}=k$ and $s_{1}^{+}(I)= \jst $: {\bf stop.} \\
    Otherwise set $s=s_{1}$ and $t=s_{1}^{+}$ and calculate $g^{(1)} =
    (f_{t}-f_{s} - 2h)/(x_{t}-x_{s})$. \\
    If $s > k$ and $g_{s^{-}(I) \, s} > g$: go to Step~6. \\
    If $t< \jst$ and $g_{t \, t^{+}(I)} > g$: go to Step~6. \\
    Otherwise: {\bf stop.}
\item Set $s_{1}=k$. Delete all elements of $I$ lying
strictly between $k$ and $\jst$ and then apply {\em
closejoin}(1). \ealg

Now define the vector $\yy(S)$ from $S$ and $h$ as follows. Let
$q'=|S|-1$ and $t_{0}=1$. When $q' \geq 1$ define $t_{ \alpha } =
s^{+}_{ \alpha } (I)$ when $h>0$ and $t_{ \alpha } = s_{ \alpha
}$ when $h=0$, for $ 1 \leq \alpha \leq q'$. Then let
\begin{equation}\label{240444}
 \Pal = [ \talm , \sal ] \qquad \mbox{ when } 1 \leq \alpha \leq
 q'+1,
\end{equation}
and define $\yy(S)$ by (\ref{240301})--(\ref{240303}). Set
$\yy^{2}$ to the value of $\yy (S)$ on exit from Algorithm
\thealg.
%
%

Step~5 is designed to avoid calculating gradients unnecessarily.
The following lemma will be used to justify this economy and also
to show that $\yy ^{2}$ is feasible. Recall the parallelogram
$\Pi(h)$ defined in the proof of Lemma~3, and given $I$ and $S$,
let $s= \sal$ and $t= s^{+}_{ \alpha } (I)$ and define
$\Pi_{\alpha }^{S} (h)$ as the parallelogram with vertices $(
x_{s}, f_{s} )$, $(x_{s}, f_{s}+ 2 (-1)^{ \alpha -1} h)$, $(x_{t},
f_{t})$, and $(x_{t}, f_{t} - 2(-1)^{ \alpha - 1} h )$. Each such
parallelogram then defines a join gradient
\begin{equation}
g^{(\alpha)}(S,h) = \frac{ f_{t} - f_{s} - 2(-1)^{ \alpha -
1}h}{x_{t} - x_{s}} . \label{24254}
\end{equation}

\begin{lemma}
Let {\em closejoin} ($\alpha$) be called, modifying $h'$, $I'$,
and $S'$ to $h$, $I$, and $S$ respectively. Suppose that there
exists $\hb \geq 0$ such that
\begin{equation}\label{24251}
  (x_{j},f_{j}) \in \Pi_{ \alpha }^{S'}( \hb )\qquad \mbox{for } \quad \salp \leq j
\leq s'^{+}_{ \alpha } (I') .
\end{equation}
and that
\begin{equation}\label{24252}
  h' \leq \hb .
\end{equation}
Then
\begin{equation}\label{24253}
  h \leq \hb .
\end{equation}
Further, if $\gb = g^{(\alpha)}(S', \hb )$, $g^{+} = g_{ s'_{
\alpha }} ( \yy (S))$ and $g^{-} = g_{ s'^{+}_{ \alpha } (I') -1
}( \yy (S))$, then
\begin{equation} \label{24255}
(-1)^ {\alpha - 1} \min (g^{+},g^{-}) \geq (-1)^ {\alpha - 1} \gb
.
\end{equation}
\end{lemma}
\subparagraph{Proof} Assume for simplicity that $\alpha$ is odd
and write $s' = \salp$ and $t' = {s'_{  \alpha }}^{+}(I')$.

First consider (\ref{24253}). The proof is trivial unless {\em
closejoin} increases $h$. In this case, $h>0$ and so there exist \jst , $k$,
$k^{+} \in I$ such that
\begin{equation} \label{24256}
    2h = \fjs - \fjs(k,k^{+} ),
\end{equation}
where $s' \leq k < \jst < k^{+} \leq t'$.
Let $\tal = s^{+}_{\alpha}(I)$.
It follows from Lemma~3
that either $\jst \leq \sal$ or $\jst \geq \tal$. The two cases
are entirely similar: it suffices to consider the first. Define
$\yb_{j}$, $s' \leq j \leq t'$, by $\yb_{s'} = f_{s'} + \hb$,
$\yb_{t'} = f_{t'} - \hb$, and $\yb_{j} = \yb_{j}(s',t')$,
$s'<j<t'$. Then the $\yb_{j}$, $s' \leq j \leq t'$, are collinear
and it follows from (\ref{24251}) that $f_{k} \geq \yb_{k} - \hb$
and $f_{k^{+}} \geq \yb_{k^{+}}- \hb$. Therefore, since
$\yb_{k}$, $\yb_{\scriptstyle j^{\star} }$, and $\yb_{k^{+}}$ are
collinear,
\begin{equation}\label{24257}
  \yb_{\scriptstyle j^{\star} } - \hb \leq \fjs (k, k^{+} ) .
\end{equation}
It also follows from (\ref{24251}) that
\begin{equation}\label{24258}
  \fjs \leq \yb_{\scriptstyle j^{\star} } +  \hb .
\end{equation}
Addition of (\ref{24257}) and (\ref{24258}) and use of
(\ref{24256}) then gives (\ref{24253}).

Now consider (\ref{24255}). The easiest case is when $\sal > s'$.
Then $g^{+}= g_{s' i}$, where $i = s'^{+}(I)$. Let the $\yb_{j}$,
$s' \leq j \leq t'$, be defined as above. Then $f_{s'}=
\yb_{s'}-\hb$ and it follows from (\ref{24251}) that $f_{i} \geq
\yb_{i} - \hb$. Then $g_{s' i}\geq g_{s'i}(\yyb)= \gb$.
Similarly, when $\tal < t'$, then $\gb \leq g^{-}$.

The next easiest case is when $s' = \sal$ and $t' = \tal$. In
this case,
\begin{equation}
g^{+} = g^{-}  = \frac{f_{t'} - f_{s'} - 2 h }
  {x_{t'} - x_{s'}},
\end{equation}
and it follows from (\ref{24253}) that this quantity is not less
than \gb, as required to establish (\ref{24255}).

It remains to resolve the two similar cases where $s'< \sal$ and
$\tal = t'$, and where $s' = \sal$ and $\tal < t'$. It suffices
to consider the former case and establish that $\gb \leq g^{-}$.
The procedure {\em closejoin} will add one or more indices
$s'^{+}$, ...,$i$, $s$ to $I$. (The gradient $g^{-}$ is now the
new join gradient, and the result required is analogous to
Lemma~4, but this lemma cannot be used directly unless $h>h'$,
which may not be the case.) The gradients defined by the
successive elements of $I$ from $s'^{+}$ to $s$ will be
monotonically increasing from $g_{s'}$, and it has already been
shown that $\gb \leq g_{s'}$, so it follows that
\beq
\gb \leq g_{is}.
\label{255a}
\eeq
Since $\tal = t'$, the test that allows Step~4 of
Algorithm~2 to adjoin $s$ to $I$ will be the inequality
\begin{equation}
  f_{t'}(i,s) \leq f_{t'} - 2h,
\end{equation}
and this inequality is equivalent to
\begin{equation}
  g_{is} \leq g,
\end{equation}
where $g$ is the new join gradient given by $(x_{t'}-x_{s})g=
f_{t'}-f_{s}-2h$. But $g=g^{-}$,
so from (\ref{255a}), $\gb \leq g$ and
(\ref{24255}) follows.~$\quad \Box$

The proof of the effectiveness of Algorithm \thealg\ and also of
its generalization in the next section will require an important
corollary of Lemma \thelemma . When new sections are added on
either side of an existing section, the convexity at the point
where they are joined increases away from zero. Thus when Steps~2
to 6 of Algorithm~4 are applied, the constraints $c_{k-1}$, $c_{
\jst -1}$, and $c_{k^{+}-1}$ increase away from zero, so that no
more than two sign changes can be created in the second divided
differences. More precisely, when defined, these constraints
satisfy
\begin{eqnarray}
  c_{k-1}( \yy ^{2}) & \geq & c_{k-1} ( \yy ^{0})  \geq  0, \label{corrlemma51}\\
c_{ \jst -1}( \yy ^{2}) & \leq & c_{ \jst -1} ( \yy ^{0})  \leq  0, \label{corrlemma52}\\
c_{k^{+}-1}( \yy ^{2}) & \geq & c_{k^{+}-1} ( \yy ^{0})  \geq  0.
\label{corrlemma53}
\end{eqnarray}

The effectiveness of Algorithm~4 can now be established.
\begin{theorem}
Algorithm~4 produces a real number $h \geq 0$, an index set $I$,
a piece set $S$ of $I$, and a vector $\yy ^{2} = \yy (S)$ well
defined by (\ref{240444}), (\ref{240301}) and (\ref{240303}) such
that $h = h(S) = F( \yy ^{2}) = \inf F(Y_{2})$ and $\yy^{2} \in
Y_{2}$.
\end{theorem}
\paragraph{Proof}
If the Algorithm stops in Step~1, then $h=h^{+}(1,n)=0$,
$I=I^{+}(1,n)$, and $S = \{ n \}$. Then $h(S)=0$. It follows from
Theorem~2 that
$\ff \in Y^{+}_{0} \subset Y_{2}$. Equation (\ref{240444}) will set $P_{1}
= [ 1,n ]$ and (\ref{240301}) and (\ref{240303}) will set
$\yy^{2} = \ff$. The theorem is then immediately established.

Otherwise the quantities \jst, $k$, and \kpl\ are well-defined and
satisfy the inequalities $1 \leq k < \jst < \kpl \leq n$. (Each
equality is possible, for example when $\ff \in Y^{\pm }_{1}$.) It
follows from Lemma~1 that at this point,
\begin{equation}\label{Red1}
  I = I^{+}(1,k) \cup I^{+} ( \kpl ,n).
\end{equation}
Step~2 is then entered. It inserts \jst\ into $I$, increases $S$ to $\hat{S} = \{ k,
\jst, n \}$, and calculates $\hat{h}= \max (h^{+}(1,k),h^{+}(
\kpl,n))$. Note that \jst, $k$, and \kpl\ are now consecutive
elements of $I$.

Steps~3 and 4 are then executed, calling {\em closejoin} in the
ranges $[k, \jst]$ and $[\jst, \kpl]$, in general increasing
$s_{1}$, $s_{2}$, and $h$, and adding new elements to $I$. The new
elements of $S$ always satisfy the inequalities $k \leq s_{1}
\leq \jst \leq s_{2} \leq \kpl$, and $h$ cannot decrease, so
$\hat{h} \leq h $. Note that $s_{2}$ can be increased to
$n$, for example when $\ff \in Y^{+}_{1}$.

Now consider the situation when the Algorithm stops. If
$h^{(1)}=0$, the Algorithm either stops in Step~4 when $h=0$, or
alternatively jumps straight from Step~5 to Step~6 re-calling {\em
closejoin}(1) with $h>0$. It follows from the properties of
Algorithm~3 that when {\em closejoin} is called with $h>0$ it
cannot increase $s$ to $t$. Therefore when the algorithm stops,
if $h>0$ then $s_{1}< \jst$ and $s_{2}< \kpl$. Now define $t_{1}$
and $t_{2}$ as for (\ref{240444}). Then $\yy^{2}$ is
well-defined in all cases. It must now be shown that it is feasible and
optimal.

The first main step is to establish (\ref{240300}), i.e.
\begin{equation}\label{Iform}
  I = I^{+}(1, s_{1}) \cup I^{-} (t_{1}, s_{2}) \cup
  I^{+}(t_{2},n).
\end{equation}
First consider the range $[1,s_{1}]$. Since $k \leq s_{1}$, $[1,
s_{1}] = [1,k] \cup [k,s_{1}]$. Let $I_{1}=I \cap [1, s_{1}]$. It
follows from (\ref{Red1}) and Theorem~3 that $I_{1} =
I^{+}(1,k) \cup I^{+}(k,s_{1})$. It is trivial that
$I_{1}=I^{+}(1,s_{1})$ unless $1<k<s_{1}$. In this case there
will exist neighbouring indices $k^{-}$ and $i$ of $k$ in $I$
such that $1 \leq k^{-} < k < i \leq s_{1} \leq \jst$. Let $h_{0}=h^{+}(1,n)$ and $g_{0} =
g^{(1)}(\hat{S},h_{0})$. Since $k^{-}$, $k$, and \kpl\ are
neighbours in $I^{+}(1,n)$, then $g_{k^{-}k} \leq g_{k \kpl} =
g_{0}$. Note that $g^{(2)}(\hat{S},h_{0}) = g_{0}$.
It is now possible to apply Lemma~5 successively with $\hb =
h_{0}$. The definition of \jst\ allows a first application of the
lemma with $h'= \hat{h}$ $S' = \hat{S}$  in the range $[k, \jst
]$ (i.e. with $\alpha=1$) to infer that at entry to Step~4,
$h^{(1)} \leq h_{0}$. This inequality and the definition of \jst\
then allow a second application of the lemma in the range $ [
\jst , \kpl]$, where also $\gb = g_{0}$, to infer that $h \leq
h_{0}$. If Step~6 is not entered, it follows from the first
application of the lemma that $g_{ki} \geq g_{0}$. If {\em
closejoin} is called again in Step~6, a third application of the
lemma may be made, in the range $[k, \jst]$, to yield that in
this case also, $g_{ki} \geq g_{0}$. Then $g_{k^{-}k} \leq
g_{ki}$. Lemma~2 can now be applied to prove that $I^{+}(1,k) \cup
I^{+} (k, s_{1}) = I^{+} (1, s_{1})$. Thus in all cases
\begin{equation}\label{Red2}
  I_{1} = I^{+}(1,s_{1}).
\end{equation}
In the same way, if
$I_{3} = I \cap [t_{2},n]$, then $I_{3} = I^{+}(t_{2},n)$, trivially when $t_{2}=n$
and otherwise by a single application of Lemma~5.

Now consider the range $[t_{1},s_{2}]$. Let $I_{2} = I \cap
[t_{1},s_{2}]$. Since $\jst \in I$, $t_{1}$ can never exceed \jst\
and so $[t_{1},s_{2}] = [t_{1}, \jst] \cup [\jst, s_{2}]$. In all
cases $I_{2}= I^{-}(t_{1}, \jst) \cup I^{-} ( \jst ,s_{2})$, and
it is trivial that $I_{2} = I^{-}(t_{1},s_{2})$ unless $t_{1} <
\jst < s_{2}$. In this case there will exist left and right
neighbours of \jst\ in I. Denote them by $j^{-}$ and $j^{+}$. Then
the successive applications of Lemma~5 made above establish that
$g_{j^{-} \jst } \geq g_{0}$ and that $g_{ \jst j^{+}} \leq
g_{0}$ so that Lemma~2 again applies to give that $I_{2} =
I^{-}(t_{1},s_{2})$. Equation (\ref{Iform}) is then established.

The feasibility of $\yy^{2}$ can now be proved. It is only
necessary to examine the four (or possibly fewer) join
constraints $c_{s_{1}-1}$, $c_{t_{1}-1}$, $c_{s_{2}-1}$, and
$c_{t_{2}-1}$ when $s_{1} < t_{1}$ and $s_{2} < t_{2}$.
First consider the last two constraints. Since $s_{2} < k^{+}$,
Lemma~4 shows that
$c_{s_{2}-1} (\yy^{2}) \leq 0$ whenever $s_{2}$ exceeds $\jst$.
When $s_{2}= \jst$, the corollary to
Lemma~5 applies to show that $c_{s_{2}-1} (\yy^{2}) \leq 0$.
Similarly, Lemma~4 when $t_{2} < \kpl$ and the corollary to
Lemma~5 when $t_{2}=\kpl$ show that $c_{t_{2}-1} (\yy^{2}) \geq
0$ whenever $t_{2} < n$.
When $t_{2}=n$ the feasibility of $\yy^{2}$ will follow automatically from the signs
of the other constraints.
Now consider the constraints
$c_{s_{1}-1}$, $c_{t_{1}-1}$. When Step~6 is entered or
$h^{(1)}=h$, the same reasoning applies to show that
$c_{s_{1}-1}( \yy^{2}) \geq 0$ when $s_{1}>1$ and $c_{t_{1}-1}(
\yy^{2}) \geq 0$ when $t_{1} > 1$. (If $\ff \in Y^{-}_{0}$, the
Algorithm will reduce $t_{1}$ to $1$.) When Step~6 is not entered
and $h> h^{(1)}$, Lemma~4 cannot be applied, but Lemma~5 and its
corollary then show that feasibility is assured unless $s_{1}>k$
or $s_{1}^{+}< \jst$. In these cases all the gradients calculated
are well-defined and the test in Step~5 allows the algorithm to stop only
when $\yy^{2}$ is feasible.

It is now necessary to show that $F(\yy^{2})=h$.

The construction of $\hat{h}$ and $h$, the inequality $\hat{h}
\leq h$, and (\ref{Iform}) show that $h= h(S)$ and $|y_{j} -
f_{j}| \leq h$ when $j \in [1,s_{1}] \cup [t_{1}, s_{2}] \cup
[t_{2},n]$, with equality when $j \in I$. For any other value of
$j$ in the range $[1,n]$, Lemma~3 and the inequality $h^{(1)}
\leq h$ show that $|y_{j} - f_{j}| < h$ whether or not Step~6 is
entered. Then $F( \yy^{2}) = h$.

If $h=0$, there is nothing further to prove. Otherwise $\sal
<\tal$, for $\alpha = 1,2$. Then it is possible to redefine $\jst
= \jst(S)$ and to let $k$, $\kpl$, and $\beta$ be uniquely redefined
from \jst . Then $\jst \in P_{\beta}$. Let $K= \{ k, \jst, \kpl,
s_{1} , t_{1} , s_{2} , t_{2} \}$. Then $K$ cannot have fewer
than five elements, even if $t_{1}=s_{2}$. Assume first that
$\beta = 1$. It follows by the same argument used in the proof of
Theorem~3 that $\civ{k \jst k^{+}}{\yy^{2}}=0$, $\civ{\jst
s_{1}t_{1}}{\yy^{2}} \geq 0$, and $\civ{s_{1}t_{1}t_{2}}{\yy^{2}}
\leq 0$. Then if $F( \vv ) \leq h$, $\civ{k \jst k^{+}}{\vv}<0$,
$\civ{j^{\star}s_{1}t_{1}}{\vv}
> 0$, and $\civ{s_{1}t_{1}t_{2}}{\vv} < 0$, and therefore the
consecutive divided differences of \vv\ with indices in the
subset $K$ must change sign twice starting with a negative sign,
so that by Theorem~1, $\vv \not\in Y_{2}$. If $\beta = 2$ or
$\beta = 3$, the argument is similar.~$\quad \Box$
%
%
%
\subsection{The general case}
The algorithm to  be described constructs a best approximation to
\ff\ over \Yq\ from $q+1$ alternately convex and concave pieces
joined where necessary by up to $q$ straight line joins. These
pieces are built up recursively from the pieces of a similar best
approximation over $Y_{q-2}$ by essentially the same method used
in the previous section when $q=2$. All the sections of the
$Y_{q-2}$ approximation remain in place except one determining
the value  of  $\inf F( Y_{q-2} )$ which is deleted and replaced
with a new piece of opposite convexity to the piece containing
this section, and two new joins. Starting with the value $h$ of a
best approximation over \Yq\ determined by the remaining
sections, the procedure {\em closejoin} is called in each join.
After this has been done, if $h$ has increased, the resulting
join constraints are checked and if necessary the calculation in
each join is repeated with the new value of $h$.
%
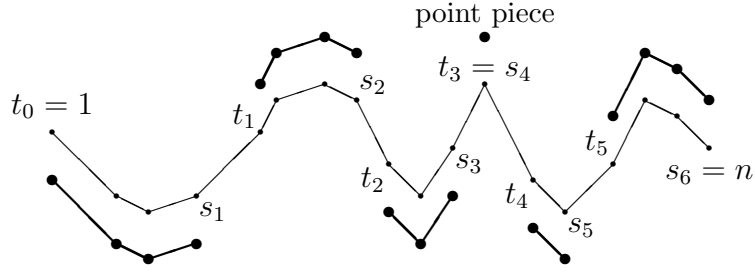
\begin{figure}
\begin{center}
\begin{picture}(258,78)
\thinlines
\dosmblob{24}{48} 
\writeat{24}{58}{$t_{0}=1$}
\doline{24}{48}{1,-1}{24} 
\dosmblob{48}{24} 
\doline{48}{24}{2,-1}{12} 
\dosmblob{60}{18} 
\doline{60}{18}{3,1}{18} 
\dosmblob{78}{24} 
\writeat{84}{18}{$s_{1}$}
\doline{78}{24}{1,1}{24} 
\dosmblob{102}{48} 
\writeat{96}{54}{$t_{1}$}
\doline{102}{48}{1,2}{6} 
\dosmblob{108}{60} 
\doline{108}{60}{3,1}{18} 
\dosmblob{126}{66} 
\doline{126}{66}{2,-1}{12} 
\dosmblob{138}{60} 
\writeat{144}{66}{$s_{2}$} 
\doline{138}{60}{1,-2}{12}
\dosmblob{150}{36} 
\writeat{144}{30}{$t_{2}$}
\doline{150}{36}{1,-1}{12} 
\dosmblob{162}{24} 
\doline{162}{24}{2,3}{12} 
\dosmblob{174}{42} 
\writeat{180}{36}{$s_{3}$} 
\doline{174}{42}{1,2}{12} 
\dosmblob{186}{66}
\writeat{186}{72}{$t_{3}=s_{4}$} 
\doline{186}{66}{1,-2}{18} 
\dosmblob{204}{30} 
\writeat{198}{24}{$t_{4}$} 
\doline{204}{30}{1,-1}{12} 
\dosmblob{216}{18} 
\writeat{222}{12}{$s_{5}$} 
\doline{216}{18}{1,1}{18} 
\dosmblob{234}{36} 
\writeat{228}{44}{$t_{5}$}
\doline{234}{36}{1,2}{12} 
\dosmblob{246}{60} 
\doline{246}{60}{2,-1}{12} 
\dosmblob{258}{54} 
\doline{258}{54}{1,-1}{12} 
\dosmblob{270}{42} 
\writeat{270}{32}{$s_{6}=n$}
\thicklines
\doblob{24}{30} 
\doline{24}{30}{1,-1}{24} 
\doblob{48}{6} 
\doline{48}{6}{2,-1}{12} 
\doblob{60}{0} 
\doline{60}{0}{3,1}{18} 
\doblob{78}{6} 
\doblob{102}{66} 
\doline{102}{66}{1,2}{6} 
\doblob{108}{78} 
\doline{108}{78}{3,1}{18} 
\doblob{126}{84} 
\doline{126}{84}{2,-1}{12} 
\doblob{138}{78} 
\doblob{150}{18} 
\doline{150}{18}{1,-1}{12} 
\doblob{162}{6} 
\doline{162}{6}{2,3}{12} 
\doblob{174}{24} 
\doblob{186}{84} 
\writeat{186}{92}{\small point piece}
\doblob{204}{12} 
\doline{204}{12}{1,-1}{12} 
\doblob{216}{0} 
\doblob{234}{54} 
\doline{234}{54}{1,2}{12} 
\doblob{246}{78} 
\doline{246}{78}{2,-1}{12} 
\doblob{258}{72} 
\doline{258}{72}{1,-1}{12} 
\doblob{270}{60} 
\end{picture}
\end{center}
\caption{Construction of $Y_{5}$ approximation}
\end{figure}
%
5
The construction of \yy\ is illustrated in Figure~6.

This method  has the disadvantage that calculations sometimes have
to be performed twice in each join, but each of the two sets of
join calculations can be performed in parallel. The procedure
{\em closejoin} is a generalization of Algorithm~2 which is a
modification of the method presented in Cullinan \& Powell
(1982). That method avoids having to repeat itself by having an
upper bound on $h$ available throughout, but because of this does
not admit of as much of its calculations being performed in
parallel. It is therefore believed that the algorithm to be
presented will often be more efficient.

Most of the notation needed for this case has already been
developed. In particular, given a piece set $S$, (\ref{240305})
defines the corresponding price $h(S)$ and when this price is
non-zero, (\ref{10601}) and (\ref{20601}) define the index \jst\
giving rise to it and the index $\beta$ of the piece within which
it lies. For any $h$ the piece \Pal\ is defined as $[\talm,
\sal]$, where $\tal = \sal^{+}(I)$ when $h>0$ and $\tal = \sal$
when $h=0$. Since every element of $I$ lies in \Pal , for some
$\alpha$, the vector $\yy(S)$ can be defined by
\begin{eqnarray}
y_{i} & = & f_{i} + (-)^{ \alpha - 1}h \qquad \mbox{ when } i \in
I \cap P_{ \alpha }, \qquad  1 \leq \alpha \leq q'+1,\label{250601} \\
 y_{j} & = & y_{j} (I) \qquad 1 \leq j \leq n. \label{250602}
\end{eqnarray}

The following algorithm then calculates $S$, $I$, and $h$ from
which $\yy^{q}$ is defined by these equations as the final value
of $\yy(S)$. Each time $h$ is increased, the new values of \jst\
and $\beta$ are calculated and recorded.
To identify joins where a second call of
{\em closejoin} will always be needed when $h$ increases from
zero, the index $\gamma$ will denote the lowest index for which
$h>0$ after the first call of {\em closejoin}($\gamma$).
When $h>0$,
the join gradients $g^{(\alpha)}(S,h)$ are defined
by (\ref{24254}).

\begin{alg} To find a best \Yq\ approximation. \end{alg}
\balg
\item Set $\qb = q$ modulo $2$. \\
      If $\qb = 0$: set $S := \{ n \}$ and use Algorithm~1 to calculate $I
         :=I^{+}(1,n)$ and $h=h^{+}(1,n)$. \\
      Otherwise set $S := \{1,n \}$, $I := \{ 1,n \}$, $h :=0$ and
        call {\em closejoin}(1).
\item If $\qb =q$ or $h=0$: {\bf stop.} \\
      Otherwise: set $j := \jst(S)$, $k := j^{-}(I)$, insert $j$
          into $I$, and insert $j$ and $k$ into $S$, increase $\qb$
          by $2$ and calculate $h := h(S)$. \\
          Set $h' := h$, $I' := I$, $S' := S$ and $\gamma := 1$.
\item  For $\alpha = 1$ to $q$:
            apply {\em closejoin}($\alpha$), and if $h$ has increased set $\alphab := \alpha$,
                                                and if $h=0$ set $\gamma := \alpha +1$. \\
      If $h=h'$: return to Step~2. \\
      Otherwise: set $\alpha = 1$ and go on to Step~4.
\item   If $\alpha < \gamma$ go to Step~5. \\
        If $\alpha = \alphab$ return to Step~2. \\
        Calculate $g = g^{(\alpha)}(S,h)$. \\
        If $\sal > \salp$ and $\mam g_{s_{\alpha}^{-}(I) \, s_{\alpha}} > \mam g$: go to Step~5. \\
        Set $t := \sal ^{+} (I')$ and $t^{+} := t^{+}(I)$. \\
        If $t < s{'_{\alpha}}^{+}(I')$ and $\mam g_{t \, t^{+}} > \mam g$: go to Step~5. \\
        Increase $\alpha$ by $1$ and repeat this step.
\item   Set $\sal := \salp$, delete all elements of $I$ lying
            strictly between \sal\ and $\sal^{+}(I')$, and
            call {\em closejoin}($\alpha$). \\
        Increase $\alpha$ by $1$ and return to Step~4.
\ealg

In practice, in order to anticipate rounding errors, the tests
whether $h=0$ should be whether $h \leq 0$. The applications of
{\em closejoin} in each step could be performed simultaneously.
In this case the largest of the ensuing values of the parameter
$h$ in Step~3 should be the value of $h$ at entry to Step~4 and \alphab\ should be set to $n$.
It will be shown in the proof of the following theorem that $h$ is
constant during Step~4.
\begin{theorem}
Algorithm \thealg\ produces a real number $h \geq 0$, an index
set $I$, a piece set $S$ of $I$, and a vector $\yy ^{q} = \yy
(S)$ well defined by (\ref{250601}) and (\ref{250602}) such that
$h = F( \yy ^{q}) = \inf F(Y_{q})$ and $\yy^{q} \in Y_{q}$.
\end{theorem}
\paragraph{Proof}

The proof that $\yy \in \Yq$ and $F(\yy )=h$ is by induction on
$\qb$. Assume that at any entry to Step~2 leading to Step~3, i.e.
with $\qb < q$ and $h>0$, that
%
\begin{equation}\label{I1}
  I =
\bigcup_{\alpha=1}^{\qb +1} I^{(-)^{\alpha}}( \talm , \sal ),
\end{equation}
where $\tal = \sal^{+}(I)$, and that
\begin{equation}\label{I2}
  h=h(S).
\end{equation}

The vector $\yy(S)$ is then well defined. Assume that
\begin{equation}\label{I3new}
  F( \yy (S))=h
\end{equation}
and that
\begin{equation}\label{I4}
  \yy(S) \in Y_{\qb}.
\end{equation}

It will be deduced from these equations that when the algorithm terminates $\yy(S) \in
Y_{q}$ and that $h = \inf Y_{q}$.
%

It has been shown in Sections~2.1 and~2.2 that \rfeqs{I1}{I3new}
hold at first entry to Step~2.

Most of the work needed for the proof has already been done in
Lemmas~3~to~5. The main task is to examine the way $h$ changes,
so as to be able to apply these lemmas. Suppose that Step~2
begins with $h=\hb$, $I= \Ib$ and $S= \Sb$, and that it increases
\qb\ to $q'$.
Step~2 also modifies $S$ from \Sb\ to $S'=\Sb \cup \{ j,k \}$
and $I$ from $\Ib$ to $I'= \Ib \cup \{ j \}$, and recalculates $h$ as $h' = h(S')$ using
(\ref{240305}). Now by (\ref{I2}) and the definition of $j$ in
Step~2, $\hb= h^{\sigma}(t_{\beta-1},s_{\beta})$, where $\beta$ is
defined in Step~2 by (\ref{20601}), $\sigma = (-)^{\beta}$, and
$t_{\beta -1} = s_{ \beta -1}^{+}(I)$.
By Lemma~1,
$I^{\sigma}(t_{\beta-1},s_{\beta})=I^{\sigma}(t_{\beta-1},k)\cup
I^{\sigma}(\kpl,s_{\beta})$, and by definition,
$h^{\sigma}(t_{\beta-1},s_{\beta}) \geq
\max(h^{\sigma}(t_{\beta-1},k), h^{\sigma}(\kpl,s_{\beta}))$. It
follows that
\begin{equation}\label{hprimehb}
  h' \leq \hb.
\end{equation}
Furthermore, it follows from the definitions of $\hb$ and $\beta$ when
$\alpha= \beta, \beta+1$ and otherwise from \rfeqs{I1}{I3new} that
\begin{equation}\label{Lemma5cdn}
  (x_{j},f_{j}) \in \Pi_{ \alpha }^{S'}( \hb )\qquad \mbox{for } \quad \salp \leq j
\leq s'^{+}_{ \alpha } (I').
\end{equation}

 Step~3 then applies {\em closejoin}($\alpha$) in each join. Let $h^{(\alpha)}$ be the
value of $h$ after {\em closejoin}($\alpha$) is called. Clearly
the $h^{(\alpha)}$ are monotonically non-decreasing.
Lemma~5 can now be applied successively beginning with (\ref{hprimehb})
to show that
\[
h^{(\alpha)} \leq \hb \qquad \mbox{ for } 1 \leq \alpha \leq q'.
\]
Let $ \talpm = s_{\alpha -1}^{+}(I')$. Then
Lemma~5 and its corollary also show that whenever $\talpm < \salp < \sal$, the
gradients on either side of $f_{\salp}$ have the correct
monotonicity for Lemma~2 to yield that $I^{\sigma} (\talpm,\salp)
\cup I^{\sigma}(\salp,\sal) = I^{\sigma} (\talpm,\sal)$, and that
when $\talm < \talpm < \salp$ that $I^{\sigma}(\talm,\talpm) \cup
I^{\sigma}(\talpm,\salp) = I^{\sigma}(\talm,\salp)$, where
$\sigma = (-)^{\alpha -1}$. It follows either from this or
otherwise trivially that in all cases after {\em
closejoin}($\alpha$) is called,
\[
I \cap \Pal = I^{\sigma}(\talm,\sal).
\]
Because $h$ cannot increase further after Step~3, this holds true
whether {\em closejoin}($\alpha$) is called once only or again in
Step~5, so that when Step~2 is next entered,
\begin{equation}\label{NewI1}
  I = \bigcup_{\alpha=1}^{q' +1} I^{(-)^{\alpha}}( \talm , \sal ).
\end{equation}

The next step is to establish that when Step~2 is next entered
with $h>0$ then $h=h(S)$. Since $\hb>0$ and $k<j< \kpl$ in
Step~2, the quantity $h(S)$ to which $h$ is set by Step~2 is well
defined by (\ref{240305}) with $\tal = s^{+}_{\alpha}(I)$. Thus
\[h'=h(S'). \]
If Step~5 recalls {\em closejoin}($\alpha$) for any $\alpha$, it will not increase
$h$ further, so $h$ always attains its final value by the end of
Step~3. It must be shown that the quantity $h(S)$ is always well
defined when Step~2 is next entered. This can only fail to be the
case if there is an $\alpha$ for which $h^{(\alpha)}=0$. In such
a case it must hold that $h'=0$ and that $h$ increased during
Step~3. Then $\alphab$ will be set to the index of
the call of {\em closejoin} in which $h$ achieved its final
positive value, and the parameter $\gamma$ will be set to the
lowest index for which $h^{(\gamma)}>0$.
Clearly, then, $\gamma \leq \alphab$ and $\alpha < \gamma \leq \alphab$, so that Step~4 will
branch to Step~5, recalling {\em closejoin}($\alpha$) with a
positive $h$ so that afterwards $ \tal= \sal^{+}(I)$. The quantity $h(S)$
will then be well defined when Step~2 is entered. It follows from
the form of Algorithm~2 and (\ref{NewI1}) that whether $h$
increases from $h'$ or not, then when Step~2 is next entered,
\begin{equation}\label{NewI2}
  h = h(S).
\end{equation}

The vector $\yy = \yy(S)$ will then be well defined at the next
entry to Step~2. The next argument will establish that
(\ref{I3new}) then holds. It follows from (\ref{NewI1}) and
(\ref{NewI2}) that it is sufficient to show that for any
indices $j$ such that $\sal < j < \tal$, $|y_{j}(S) - f_{j}| <
h$. This is equivalent to the statement that the graph point
$(x_{j},y_{j}(S))$ lies inside the parallelogram
$\Pi^{S}_{\alpha}(h)$. It follows from Lemma~3 that immediately
after {\em closejoin}($\alpha$) is last applied,
$(x_{j},y_{j}(S))$ lies strictly within $\Pi^{S}_{\alpha}(h)$. If
this call of closejoin comes in Step~5, $h$ has already attained
its value at next entry to Step~2 and the result is immediately
established. If, on the other hand, this call of {\em closejoin} comes
in Step~3, then $(x_{j},y_{j}(S))$ lies inside
$\Pi^{S}_{\alpha}(h^{(\alpha)})$, and $h^{(\alpha)} \leq h$. A
simple calculation shows that if $a < b$ then
$\Pi^{S}_{\alpha}(a) \subset \Pi^{S}_{\alpha}(b)$. The result
then follows.
Thus at next entry to Step~2,
\begin{equation}\label{NewI3}
  F( \yy (S))=h.
\end{equation}

The next argument will re-establish (\ref{I4}) at that point. It
follows from (\ref{I1}) that it is only necessary to establish
that the join constraints $c_{\salm}$ and $c_{\talm}$ of \yy\
have the correct signs when $2 \leq \sal < \tal \leq n-1$, i.e. when $h>0$.

When $1<\sal=\salp$, (\ref{I4}) and the corollary to Lemma~5
imply that $c_{\salp-1}(\yy)$ had the correct sign at the last
exit from Step~2 and has not now moved closer to zero. The same is
true of $c_{\talm}$ when $\tal = \talp < n$
When $\sal > \salp$, Lemma~4 applies provided that $h$ has not
increased after the last call of {\em closejoin}($\alpha$) to
yield that $(-1)^{\alpha - 1}c_{\sal - 1} \geq 0$ and that when
$1 < \tal < \talp$ that $(-1)^{\alpha - 1}c_{\tal - 1} \leq 0$.
If {\em closejoin}($\alpha$) is only called once and $h$
increases subsequently, it must be shown that the test in Step~4
is adequate to ensure feasibility. In this case Step~4 will
certainly be entered (because $h>h'$) and Step~3 will set $\alphab$
and $\gamma$ so that $\gamma \leq \alpha < \alphab$, so that
$h^{(\alpha)}>0$. It follows that $\sal < \tal$ and so
$g$ is well defined and the correct sign of $c_{\salm}$
is assured by the test in Step~5. The case $\tal < \talp$ is
similar.
Thus at next entry to Step~2, even when $h=0$,
\begin{equation}\label{NewI4}
  \yy \in Y_{q'}.
\end{equation}

Thus \rfeqs{I1}{I4} are established by induction. It follows that
when the algorithm terminates with a number $h$ and the pieces
from which the vector \yy\ is constructed, $\yy \in Y_{ \qb } \subset \Yq$ and $F(
\yy)=h$.

\vspace{1ex}
The proof of optimality is similar to that in Theorem~4.
If $h=0$ there is nothing to prove.
Otherwise, after the algorithm has terminated,
let $\sal$, $\tal$, $1 \leq \alpha \leq q+2$,
be defined as the join points that
would next have been constructed in Step~2
if the algorithm had not terminated
(i.e. by adding $k$, \jst, and \kpl\
to the existing set of join points and reindexing.)
It follows as in the proof of Theorem~3 that because $h>0$,
$\sal < \tal$ for all $\alpha$.
Now define $K= \{ \, \sal , \, \tal \, : \,  1 \leq \alpha \leq q+2 \, \}$.
Because the piece
$[ t_{\beta} , s_{\beta+1}]$ has only one element and
any of the other pieces can have only one element,
this set $K$ can have from $q+3$ to $2q+3$ elements.
It has to be shown that the consecutive divided differences with indices in $K$
of any vector $\vv \in \realsn$ giving
a lower value of $F$ than $h$ change sign $q$ times starting
with a negative sign, so that by Theorem~1, $\vv \not\in \Yq$.

It follows from the construction of $\yy$ and (\ref{I4}) that
$ \mam c_{i-1} \yy \leq 0$ for $\sal \leq i \leq \talpl $ and so, from the corollary to Theorem~1, that
for any $\alpha$ in the range $1 \leq \alpha \leq q+2$, $\csts{\yy}  \leq  0$ and $\csst{\yy} \leq 0$.
If $\tal = s_{\alpha +1}$, it follows immediately from (\ref{250601}) that if $F(\vv)<h$
\[
\cstt{\vv} < 0.
\]
Otherwise when $\tal < s_{\alpha +1}$, it follows from (\ref{250601}) that if $F(\vv)<h$, then
$\cstt{\vv} < 0$ and $\csst{\vv} < 0$ and hence, from the corollary to Theorem~1, that the
{\em consecutive} divided differences with indices in $K$ satisfy
\[
\csts{\vv} < 0 \qquad \mbox{\em or} \qquad \ctst{\vv} < 0.
\]

Therefore the divided differences of \vv\ with indices in $K$ have
at least $q$ sign changes starting with a negative sign,
and so if \vv\ is any vector  for which $F(\vv) < h$,
then $\vv \not\in \Yq$.~$\quad \Box$

%
%
%
%
\section{Conclusions}
This section will describe the results of some tests of the data
smoothing method developed in the previous section.
These tests were conducted by contaminating sets of values of a known function
with errors, applying the method, and then comparing the results
with the exact function values.
If \ggg\ is the vector of exact function values, a simple measure of the
effectiveness of the method can be obtained from the quantity
\[
P_{p} = \left ( 1 - \frac{\|\yy - \ggg \|_{p}}{\|\ff - \ggg \|_{p}} \right ),
\]
obtained from the $\ell_{p}$ norm.

Algorithm \thealg\ was trivially extended to find a best approximation over
$Y_{\pm q}$.
It was coded in PASCAL and run on a SUN computer.
The errors added to exact function values were either truncation or
rounding errors, or uniformly distributed random errors in the
interval $[ - \epsilon , \epsilon ]$.
For each set of data the values of \Pu\ and \Pe\ were calculated.
Many of the results were also graphed.

As a preliminary test, equally spaced values of the zero function
on $[ -5 , 5 ]$ were contaminated with uniformly distributed errors
with $\epsilon = 0.1$. The results are shown in Table~1.
\begin{table}[hbp]
\begin{center}
\begin{tabular}{||c|c|c||} \hline
$n$ & \Pu & \Pe \\ \hline
5001 & -19.37 & 94.51 \\ \hline
501  & -68.97 & 76.25 \\ \hline
\end{tabular}
\caption{The zero function}
\end{center}
\end{table}

The difference between $y_{j}$ and zero was of the order of $10^{-4}$
for most of the range, but near the ends it rose to $10^{-1}$,
accounting for the relatively high value $\| \yy - \ggg \|_{2}=0.3098$
when $n=501$.
High end errors were frequent and so the value of \Pu\ was not usually a reliable indicator of the efficiency of
the method.

Two main types of data were then considered: undulating data and peaky data,
because it is these types of data that can be hard to smooth using divided
differences unless sign changes are allowed.
The first main category of
data were obtained from equally spaced values of the function
$\sin \pi x$ on the interval $[ -2 , 2 ]$,
and the second from equally spaced values of the normal distribution function
$\displaystyle N_{s}(x) = (2 \pi x )^{- \frac{1}{2}} \exp (-x^{2}/2s^{2})$
with $s = 0.8$, on the interval $[ -5 , 5 ]$.
These same functions were previously used to test the $\ell_{2}$
data smoothing method of Cullinan(1990) which did not allow sign changes, where it was shown that the sine data were possible to treat well but the peaked data did not give very good results.

Many of the results looked very acceptable when graphed, even when
\Pe\ was only moderate.
The results of Table~1 have already indicated that it is quite difficult for \Pe\ to
approach $100$. This is often because of end errors, as there, but there is also
another phenomenon that can move smoothed points further away from the
underlying function values.
The method raises convex pieces and lowers concave pieces, and this often
results in points near an extremum of the underlying function being
pulled away from it, for example if there is a large error on the low side of a
maximum.

Different sets of random errors with the same $\eps$ can give very
different values of both \Pu\ and $\Pe$, particularly
when the spacing between points is not very small, and in fact it was found that
reducing the spacing beyond a certain point can make a great difference.

The method coped quite well with the sine data and
was very good with the peaked data, and a great improvement on the method in
Cullinan(1990), managing to model both the peak and the flat tail well.
The results of one run are shown in Figure~7.

\begin{figure} [htbp]
\begin{center}
\resizebox{\textwidth}{!}{\includegraphics{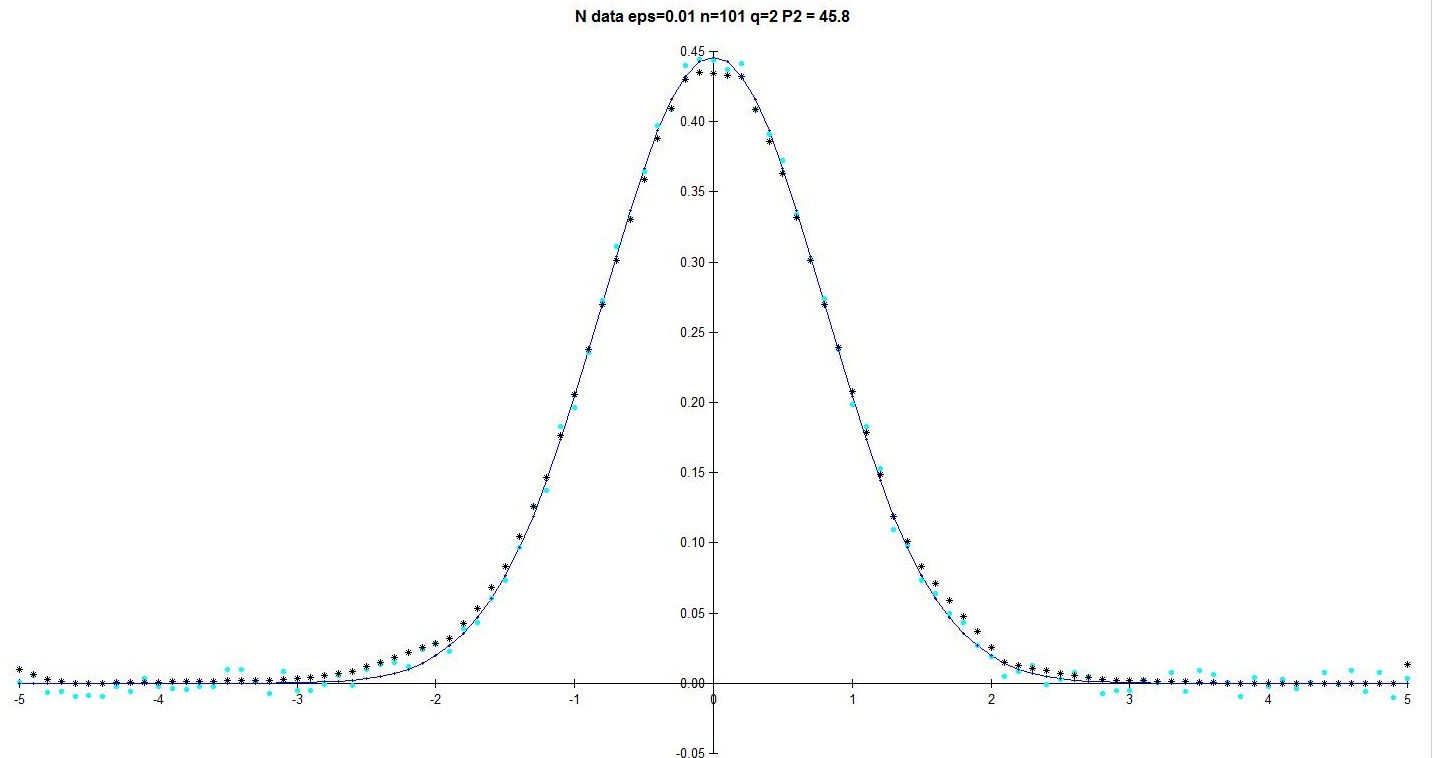}}
\end{center}
\caption{Peaked Data with $n=1$, $q=2$, $\eps = 0.01$, $\Pe = 45.8$}
\end{figure}

The method is very economical, and can give very good results, particularly
when the data are close together. Thus it seems a good candidate for large
densely packed data, even when the errors are relatively large.
BIBLIOGRAPHY
%

\listoffigures
\end{document}